\documentclass{svjour3}
\smartqed

\usepackage{diagbox}
\usepackage{epic,eepic}
\usepackage{multirow}
\usepackage{tikz}
\usepackage{amsfonts,amssymb}
\usepackage{caption}
\usepackage{bm}
\usepackage[bf,SL,BF]{subfigure}
\usepackage{epstopdf}
\usepackage{graphicx,ifthen,latexsym,mathrsfs,multicol}
\usepackage{amsmath}
\usepackage{epstopdf}
\usepackage{bbding}
\usepackage{mathptmx}
\usepackage{algorithm}
\usepackage{algorithmic}
\usepackage{color}
\usepackage{subfigure}
\usepackage{appendix}
\usepackage{booktabs}
\usepackage{algorithm}
\usepackage{algorithmic}
\usepackage{mathtools}
\newcommand{\tabincell}[2]{\begin{tabular}{@{}#1@{}}#2\end{tabular}}
\numberwithin{equation}{section}
\numberwithin{theorem}{section}
\numberwithin{lemma}{section}
\numberwithin{remark}{section}

\begin{document}
\title{Discovery of subdiffusion problem with  noisy data via deep learning}

\author{Xingjian Xu    \and
         Minghua Chen
}

\institute{X. Xu \and M. Chen (\Envelope)   \at
              School of Mathematics and Statistics, Gansu Key Laboratory of Applied Mathematics and Complex Systems, Lanzhou University, Lanzhou 730000, P.R. China\\
email:chenmh@lzu.edu.cn;
}

\date{Received: date / Accepted: date}

\maketitle

\begin{abstract}
Data-driven discovery of partial differential equations (PDEs)  from observed data in machine learning  has been developed by embedding the discovery problem.
Recently, the discovery of traditional ODEs dynamics using linear multistep methods in deep learning have been discussed in
[Racheal and Du, SIAM J. Numer. Anal. 59 (2021)  429-455; Du et al. arXiv:2103.11488].
We extend this framework to the data-driven discovery of the   time-fractional PDEs, which can effectively characterize the ubiquitous power-law phenomena.
In this paper, identifying source function of subdiffusion  with noisy data using  $L_{1}$ approximation in deep  neural network is presented.
In particular, two types of networks for  improving the generalization of the  subdiffusion problem are designed with noisy data.
The numerical experiments  are given to illustrate  the availability using deep learning.
To the best of our knowledge, this is the first topic on the discovery of  subdiffusion in deep learning with noisy data.
\end{abstract}

\section{Introduction}
Deep learning has been extended to many different practical fields, including image analysis, natural language processing, system recognizing etc \cite{IYA:2016}.
There are already some important progress for numerically solving the ODEs and PDEs systems in deep learning.
For example,  variational methods \cite{WEBY:2017,YGMK:2021,YLJ:2021,CGD:2021}, Galerkin methods \cite{DGM:2018,CJR:2021}, random particle method \cite{XUY:2020}, different operators method \cite{LZY:2021}, linear multistep method \cite{DUQ:2021}. Based on the the neural networks, the data-driven discovery of differential equations has been proposed in \cite{RTP:2019,SHR:2019,DEM:2020,MRI:2018,QWX:2019,DUQ:2021} including the space-fractional differential equations \cite{GRPK:19}.
In particular,  the discovery of traditional ODEs dynamics using linear multistep methods in deep learning have been discussed in \cite{DUQ:2021,RD:21} and  \cite{RTP:2019}.
Identifying source function from observed data is a meaningful and challenge topic for the time-dependent PDEs.

Over the  few decades, fractional models have been attracted wide interest since it can effectively characterize the ubiquitous power-law phenomena
with noise data \cite{EK:11,LWD:17},
which are applied in underground environment problems, transport in turbulent plasma, bacterial motion transport, etc.  \cite{M2000}.
In this work, we study the discovery of  the following subdiffusion  in deep learning with noise data, whose  prototype is \cite{LWD:17,MS:12,WCDBD:22}, for $0<\alpha \leq  1,$
\begin{equation}\label{ad1.1}
\begin{cases}
{_{0}^CD_t^{\alpha}u(x,t)}+A u(x,t)=f(x,t):=\overline{f}(x,t)+ {\rm perturbation ~(noise) }, \\
u(0,t)=0,~~~u(l,t)=0,~~~ t\in (0, T], \\
u(x,0)=g(x),~~~x\in [0, l]\\
\end{cases}
\end{equation}
with  $A$ a positive definite, selfadjoint, linear Laplacian operator. The Caputo fractional derivative \cite{Podlubny:99} is defined by
\begin{equation}\label{1.2}
 _0^{C}\!D_{t}^{\gamma}u(x,t)=\left\{ \begin{array}
 {l@{\quad} l}
\displaystyle \frac{1}{\Gamma(1-\alpha)} \int_{0}^t \frac{\partial u(x,{\eta})}{\partial {\eta}} {(t-\eta)^{-\alpha}}d\eta,& 0<\alpha <1,\\
 \cr\noalign{\vskip 0 mm}
    \displaystyle \frac{\partial u(x,t)}{\partial t},&\alpha =1,
 \end{array}
 \right.
\end{equation}
and $f$ is a given forcing term $\overline{f}$ with noise data.
In practice, $f$ is not completely known   and is disturbed around some known quantity $\overline{f}$, i.e.,
$$f(x,t)=\overline{f}(x,t)+ {\rm perturbation ~(noise) }.$$
Here $\overline{f}$ is deterministic and is usually known while no exact behavior exists of the perturbation (noise) term.
The  uncertainty (lack of information) about $f$ (the perturbation term) is naturally represented as a stochastic quantity, see  \cite{ZK:17}.
In the following, we focus on the uniform noise and Gaussian noise. More general noise  \cite{ZK:17} such as white noise,
Wiener process or Brownian motion (which corresponds to the stochastic fractional  PDEs \cite{LWD:17,WCDBD:22})  can be similarly studied,
since it can be discretized  as like uniform noise form  by a simple arithmetic operations \cite{LWD:17,MS:12,WCDBD:22}.

The main contribution  of this paper is to discover  subdiffusion via $L_{1}$ approximation in deep learning with  noisy data,
where we design two types  neural networks to deal with the fractional order $\alpha$, since $\alpha$ is a variable parameter.
More concretely, we use two strategies  (fixed  $\alpha$ and variable $\alpha$, respectively)   to train network, which improves  the generalization for the  subdiffusion problem.
The advantage of  first Type  is that model \eqref{C6} can be  trained    for each fixed $\alpha$, which reduces the computational count and required storage, since it fixes  $\alpha$ as input data.
Obviously, it may loss some accuracy. To recover the accuracy,  second Type  is complemented, which  needs more computational count.
It implies  an interesting generalized structure by combining  Type 1 and Type 2, which may keep suitable accuracy and reduces computational count.

The paper is organized as follows. In the next section, we introduce  the discretization schemes for the subdiffusion problem  \eqref{ad1.1} and its discovery of subdiffusion in Section 3.
We construct the  discovery of subdiffusion based on the basic deep neural network (DNN)  in Section 4. The two types DNN are  designed/developed for the subdiffusion model in Section 5.
To show the effectiveness of the presented schemes, results of numerical experiments with   noisy data   are reported in Section 6.
Finally, we conclude the paper with some remarks on the presented results.
	
\section{Subdiffusion problem}
In this section, we introduce the  $L_{1}$ approximation for solving the subdiffusion.
Let $t_{n} = nh_{t}$, $n=0,1,...N_{t}$ be a partition of the time interval $[0,    T]$ with the grid size $\tau  = T/N_t$. Denote  the mesh points $x_m= m h_{x}$, $m=0, 1, \cdots, N_{x}$ with the uniform grid size $h_{x} =l/M$, $\Omega=(0,l)$.

The objective for solving the initial value problem in \eqref{ad1.1} is to find the approximation $u_{m}^{n} \approx u(x_{m},t_{n})$ for given source function $f(x_{m},t_{n})$.

The diffusion term is approximated by a standard second-order discretization:
\begin{equation}\label{A1}
u_{xx}(x_m, t_{n})\approx \delta_x^2 u_m^{n}:=\frac{u_{m+1}^n-2u_{m}^n+u_{m-1}^n}{h^2}.
\end{equation}

There are serval ways to discretize the  Caputo fractional substantial derivative \eqref{1.2}.
For convenience,  we use the following standard  $L_1$ scheme \cite{ChJB:21,SOG:17,Lin:07} with the truncation error  $\mathcal{O} \left( \tau^{2-\alpha} \right)$, namely,
\begin{equation}\label{A12}
\begin{split}
_{0}^CD_t^{\alpha} u(x_m, t_n)\approx D_\tau^{\alpha} u_m^n
=\frac{1}{\Gamma(1-\alpha)}\!\sum\limits_{k=0}^{n-1}\frac{u_m^{k+1}\!-u_m^k}{\tau }
 \int_{t_k}^{t_{k+1}}\!(t_n-s)^{-\alpha}ds
 = \sum_{k=0}^{n}\omega^{(\alpha)}_{n-k}u_{m}^{k}.
\end{split}
\end{equation}
Here the  coefficients are computed by $\omega^{(\alpha)}_{0} = \frac{1}{\Gamma(2-\alpha)\tau ^{\alpha}}$,
$\omega^{(\alpha)}_{n} =\frac{(n-1)^{1-\alpha}-n^{1-\alpha}}{\Gamma(2-\alpha)\tau ^{\alpha}}$ and
\begin{equation*}
\begin{split}
\omega^{(\alpha)}_{n-k}=\frac{(k+1)^{1-\alpha}-2k^{1-\alpha}+(k-1)^{1-\alpha}}{\Gamma(2-\alpha)\tau ^{\alpha}},~~ 1 \leq k \leq n-1.
\end{split}
\end{equation*}
Thus, we approximate \eqref{ad1.1} by the discrete problem with given $f(x_{m},t_{n}) = f_{m}^{n}$
\begin{equation}\label{AA3}
\begin{split}
\sum_{k=0}^{n}\omega^{(\alpha)}_{n-k}u_{m}^{k}-\delta_x^2 u_m^n=f_{m}^{n}.
\end{split}
\end{equation}

\begin{remark}
When $0<\alpha<1$,  the time Caputo fractional derivative  uses the information of the classical derivatives
 at all previous time levels (non-Markovian process).  If  $\alpha=1$,  it can be seen that by taking the limit $\alpha \rightarrow 1$ in (\ref{1.2}), which gives the following equation
\begin{equation*}
\frac{\partial u(x,t_{n})}{\partial t} =\frac{u(x,t_{n})-u(x,t_{n-1})}{\tau} + \mathcal{O}(\tau).
\end{equation*}
\end{remark}
\section{Discovery of subdiffusion}
The discovery of subdiffusion is essentially an inverse process of solving a subdiffusion problem (\ref{1.2}).
That means, suppose that only the information of $u(x_{m},t_{n})$ at the pairs of the uniform grid points $\left\{x_{m},t_{n}\right\}$ is provided, we need to recover the source function $f$ \cite{RD:21}.
Assume $u(x,t)$ and $f(x,t)$ both are unknown in subfiffuion system \eqref{ad1.1} with given $u_{m}^{n}=u(x_{m},t_{n})$, the target is to approximate the close-form expression for $f(x_{m},t_{n})$.
We can use the approximated source function $f_{m}^{n}$ to rebuild the discrete relation between $f_{m}^{n}$ and $u_{m}^{n}$, namely,
\begin{equation}\label{CC3}
\begin{split}
f_m^n=\sum_{k=0}^{n}\omega^{(\alpha)}_{n-k}u_{m}^{k}-\delta_x^2 u_m^n.
\end{split}
\end{equation}

It should be noted that \eqref{CC3} directly follows the $L_1$ approximation in (\ref{AA3}). Different from \eqref{AA3} that evaluates $u_m^n$ given $f_m^n=f(x_m,t_n)$,
model \eqref{CC3} computes $f_m^n$ from the given data $u_m^n$. It also indicates the discovery of subdiffusion is actually an inverse process of solving the subdiffusion problem.

\section{Neural Network approximation}
In this section, we  introduce the basic deep neural network (DNN) \cite{CWHZ:21,ZuowShen:2021,CGD:2021,SHR:2019}, which will be used in discovery of subdiffusion.
\subsection{Structure of deep neural network}\label{sub4.1}
As the input vector $\bm{x}$, the neural network can be denoted as
\begin{equation*}\label{C4}
\hat{\phi} (\bm{x};\bm{\theta})= L_{M} \circ \sigma \circ L_{M-1}\circ\sigma ...\sigma \circ L_{1}(\bm{x})
\end{equation*}
\begin{figure*}[htbp]
	\centering
	\includegraphics[scale=0.40]{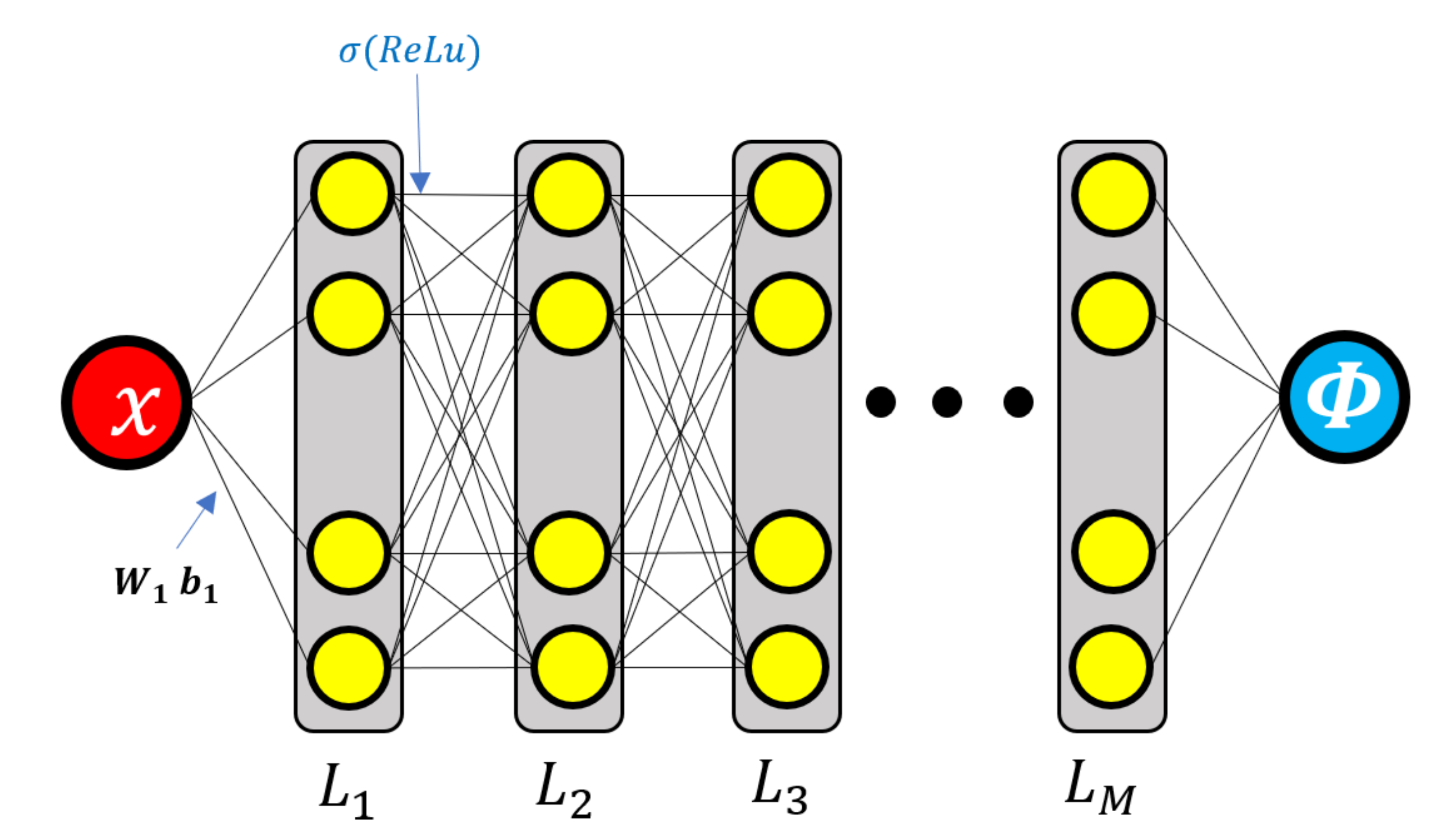}
	\caption{\textbf{Structure of deep neural network (DNN)}}\label{Fig.1}
\end{figure*}
where $\hat{\phi}$ is a nonlinear function in the whole neural network. The integer $M$ is number of the layer of DNN. The $M$-th hidden layer has the following affine transform:
\begin{equation}\label{ad4.1}
	L_{M}(z) = \bm{W}_{M}z+\bm{b}_{M}, 1 \le m \le M
\end{equation}
where $\bm{W}_{M},\bm{b}_{M}$ belongs to a parameter functional space $\theta$. The weight family $\bm{W}_{M} \in \mathbb{R}^{p_{m}\times p_{m-1}}$ and collection of bias $\bm{b}_{M} \in \mathbb{R}^{p_{m}}$.
The dimensional $p_{m} $ is the number of neurons (width of neural network) in the m-th layer.  The activation function $\sigma$ is used  as non-linear processing of information. Here we choose the ReLU function $\sigma(t)=\max\left\{0,t\right\}$ as  the activation function. All DNN structures can be presented in Figure \ref{Fig.1}.

\subsection{Discovery of subdiffusion  in deep learning}

Among many different structures of approximations, it is more predominant to employ neural networks as a nonlinear approximation tool to get source function $f(x,t)$, which is convenient to compute and implement. So the neural network approximation is focused in this paper. Consider the neural network approximation via $L_{1}$ discertization. We use $\mathscr{N}$ as the set of all neural networks with special architecture  \cite{RD:21}. Now we introduce a network $\widetilde{f}(\cdot) \in \mathscr{N}$ to approximate $f(\cdot)$.
The model  of  neural network approximation is developed by replacing each $f_{m}^{n}$ with  $\widetilde{f}_{m}^{n}$ in \eqref{CC3}
\begin{equation} \label{C6}
\widetilde{f}_{m}^{n} =\sum_{k=0}^{n}\omega^{(\alpha)}_{n-k}u_{m}^{k}-\delta_x^2 u_{m}^{n}.
\end{equation}

Now we seek $\widetilde{f}(\cdot)$ by minimizing  the residual error \eqref{C6} under a machine learning framework in practice,  namely,
\begin{equation}
J_{h}(\widetilde{{f}}) =\min_{u \in \mathscr{N}} J_{h}(u)
\end{equation}
with  the loss function
\begin{equation}\label{ad4.3}
J_{h}(u) =\frac{\sum_{m=1}^{N_{x}}\sum_{n=1}^{N_{t}}|\widetilde{f}_{m}^{n}- \sum_{k=0}^{n}\omega^{(\alpha)}_{n-k}u_{m}^{k}+\delta_{x}^{2}u_{m}^{n}|^2}{N_{t}N_{x}-1}.
\end{equation}

\section{Two types of  deep neural networks}
It is an interesting question  how we can train the discovery of subdiffusion  with fractional order $\alpha$, since  $\alpha$ is variable parameter.
Developed  the structure of deep neural network in Subsection \ref{sub4.1},
we  use two types   (fixed  $\alpha$ and variable $\alpha$, respectively)  of DNN for training the discovery of subdiffusion  \eqref{C6} in this work.

\subsection{Type 1 of DNN: fixed  $\alpha$}
\begin{figure}
\centering
\includegraphics[scale=0.35]{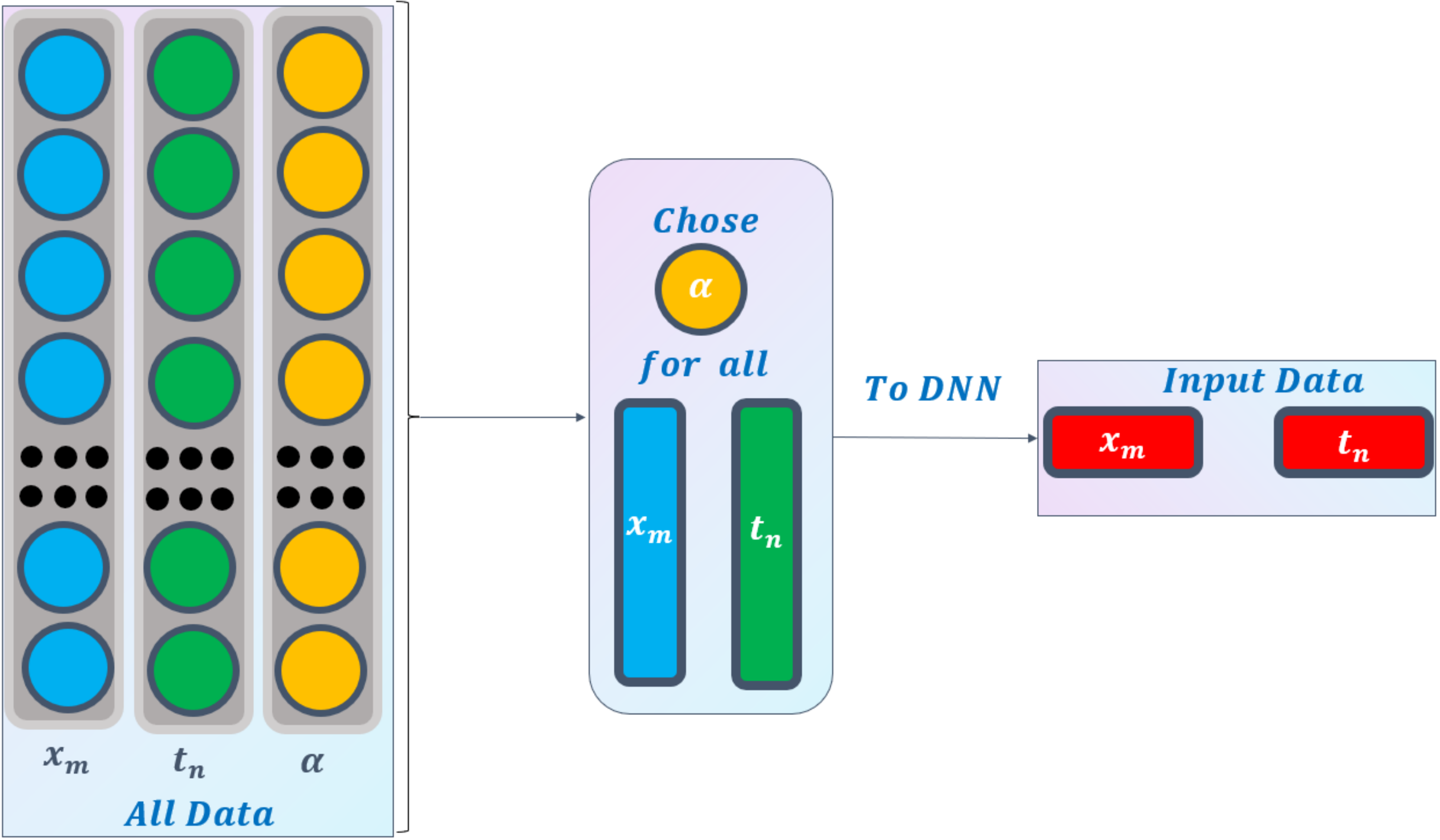}
\caption{\textbf{The construction for Type 1 Network.}}\label{Fig.ad2}
\end{figure}

To  train   the source function $\widetilde{f}$ of subdiffusion model \eqref{C6} by the deep neural network approximation,
we design the Type 1 of DNN for the fixed $\alpha$, see Fig \ref{Fig.ad2}, which just needs $(x_m,t_n)$ as input data.
\subsection{Type 2 of DNN: variable   $\alpha$}
\begin{figure}[t]
\centering
\includegraphics[scale=0.35]{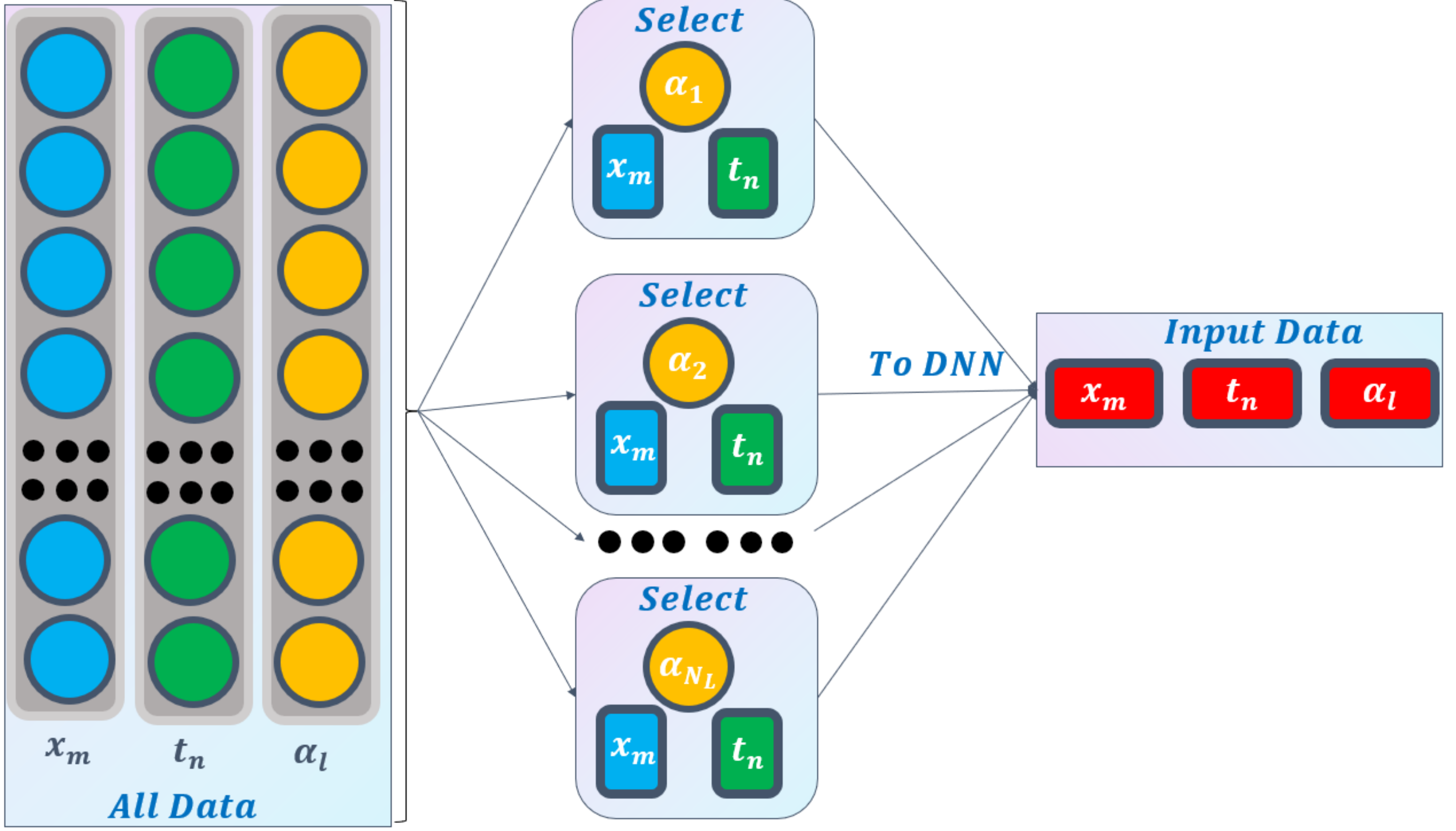}
\caption{\textbf{The construction for Type 2 Network.}}\label{Fig.ad3}
\end{figure}

To  discover/recover  the source function $\widetilde{f}$ of subdiffusion model \eqref{C6} in DNN,
we also design the Type 2 of DNN for the variable $\alpha$, see Fig \ref{Fig.ad3},  which uses  $(x_m,t_n,\alpha_l)$ as input data.
Here $\alpha$ can be chosen a sequence $\left\{\alpha_{l}\right\}_{l=1}^{N_l}$.

\begin{remark}
The advantage of  Type 1 is that model \eqref{C6} can be  trained    for each fixed $\alpha$, which reduces the computational count and required storage, since it fixes  $\alpha$ as input data.
Obviously, it may loss some accuracy.
To recover the accuracy,  Type 2 is complemented, which  needs more computational count.
Hence, there is an interesting generalized structure by combining  Type 1 and Type 2, which keeps suitable accuracy and reduces computational count.
\end{remark}

\section{Numerical Experiments}
In this section,  we discover subdiffusion problem \eqref{ad1.1} with uniform noise and Gaussian noise in deep learning.
 More general noise \cite{ZK:17} such as white noise,  Wiener process or Brownian motion
(which corresponds to the stochastic fractional  PDEs \cite{LWD:17,WCDBD:22})  can be similarly studied,
since it can be discretized  as like uniform noise form  by a simple arithmetic operations \cite{LWD:17,MS:12,WCDBD:22}.
Several examples are provided to show the performance of subdiffusion discovery via $L_1$ approximation in DNN by Algorithm \ref{algorithm1}.

\begin{table*}[h!]
\begin{center}
 \caption*{\textbf{The overall setting in all experiments is summarized  as following} }\label{tabad.1}
\resizebox{\textwidth}{!}{
\begin{tabular}{l|c}
	\hline
	\textbf{Feature}&\textbf{Statement}  \\
	\hline
		Environment & \tabincell{c}{All numerical experiments are programmed in Python 3.8,\\ and Pytorch library for neural network are implemented. } \\
	\hline
	\tabincell{c}{Optimizer} & \tabincell{c}{Apply Adam optimization  algorithm for solving  minimization \\ problem   \eqref{ad4.3}
with the learning rate  $1.0\times 10^{-2}$.}\\
\hline
\tabincell{c}{  Parameters} &\tabincell{c}{ By fixing the epoch, the batch size is changed to derive the optimization.}\\
	\hline
\tabincell{c}{Network setting} &\tabincell{c}{100 neurons and 10 layers with  ReLU activation function in Fig. \ref{Fig.1}, \\
and uniform distribution $\bm{W}_{M},\bm{b}_{M}$ in \eqref{ad4.1}  is initialized.
}\\
\hline
\tabincell{c}{Two types  setting} & \tabincell{c}{We  fix the fractional order $\alpha=0.3$ in Type 1, \\ and choose a sequence $\alpha=\left\{l/10\right\}_{l=1}^{10}$ in  Type 2.}\\
\hline
\end{tabular}}
\end{center}
\end{table*}
The relative error $||\cdot||_r$ is measured  by
\begin{equation} \label{D6.1}
||e_{\widetilde{f} }|| = \frac{\left(\sum_{m=1}^{N_{x}}\sum_{n=1}^{N_{t}}|\widetilde{f}_{m}^{n}-f(x_{m},t_{k})|^2\right)^{1/2}}
{\left(\sum_{m=1}^{N_{x}}\sum_{n=1}^{N_{t}}|f(x_{m},t_{k})|^2\right)^{1/2}},
\end{equation}
where  we  actually   use the well-known  Frobenius norm of the $n\times n$ matrix $A$, i.e.,
$$||A||_F=\left(\sum_{i,j=1}^n|a_{i,j}|^2\right)^{1/2}.$$

 \begin{algorithm}
 \caption{Discovery of subdiffusion in deep learning}\label{algorithm1}
 \begin{algorithmic}[1]
 \REQUIRE
         Batch size $B_1$; Epoch $I_1$; Learning rate $r$; Size of DNN.
 \ENSURE Discovery $\widetilde{f}_{m}^{n}$.
 \STATE Generate mesh data $x_{m}$ and $t_{n}$ with $h_{x}=h_{t}=1/100$ on the  interval $[0,1]$.
 \STATE Calculate $u_{m}^{n}$ with $x_m$, $t_n$ and $\alpha$ by (\ref{AA3}).
 \STATE Discover $f_{m}^{n}$ by using $u_{m}^{n},x_{m},t_{n}$  in (\ref{CC3}).
 \STATE Construct the  training set $\left\{x_{m},t_{n}\right\}$ in Type 1 or $\left\{x_{m},t_{n},\alpha\right\}$ in Type 2 to approximate $f_{m}^{n}$.
 \STATE Initialize the weights $\bm{W}$ and biases $\bm{b}$ in \eqref{ad4.1} with Normal distribution of DNN,  see Fig \ref{ad4.1}.
\FOR{iteration $i=1:I_{1}$ in epoch}
\FOR{iteration $j=1:B_1$ in Batch size}
\STATE Construct DNN for training set in small batch and get the output $\widetilde{f}_m^{n}$ in \eqref{C6}
\STATE Update $\widetilde{f}_{m}^{n}$ the  neural network by minimizing the loss function (\ref{ad4.3})
\STATE Update the neural network by  minimizing the loss function in Fig \ref{ad4.1}:
\STATE $\bm{W}_{j+1} \gets  \bm{W}_{j} + r \times \bm{W}_{j}$
\STATE $\bm{b}_{j+1} \gets  \bm{b}_{j} + r \times \bm{b}_{j}$
\ENDFOR
\ENDFOR
\STATE Calculate the relative error $||\cdot||_{r}$ in (\ref{D6.1}).
 \end{algorithmic}

 \end{algorithm}
\subsection{Numerical Experiments for discovery with uniform  noise}
\begin{example}\label{ex2}
Let us consider the following subdiffusion problem \eqref{ad1.1}
\begin{equation*}
	\begin{cases}
		{_{0}^CD_t^{\alpha}u(x,t)}-\Delta u =f, \quad x \in [0,1], t \in [0,1];\\
		u(x,0)=\sqrt{x(1-x)};\\
		u(0,t)=0, u(1,t)=0.		
	\end{cases}
\end{equation*}
Here the source function is given by
$$f(x,t) = (t+1)^{2}\left(1+\chi_{(0,1/2)}(x)\right)+ \delta$$
with
\begin{equation*}
    \chi_{(0,1/2)}(x) =\begin{cases}
        1,\quad 0\le x \le 1/2,\\
        0,\quad {\rm otherwise}.
    \end{cases}
\end{equation*}
Here $\delta$ is the  uniform  noise with $0\%$ (clean data) $10\%$, $20\%$, $50\%$ level for a single uniform distributed random number in the interval $(0,1)$, respectively.
In this experiment, we train  the deep learning discovery $\widetilde{f}$ in \eqref{C6}.
\end{example}

\subsubsection{Type 1 in DNN}
\begin{figure*}[!h]
	\centering
	\includegraphics[scale=0.35]{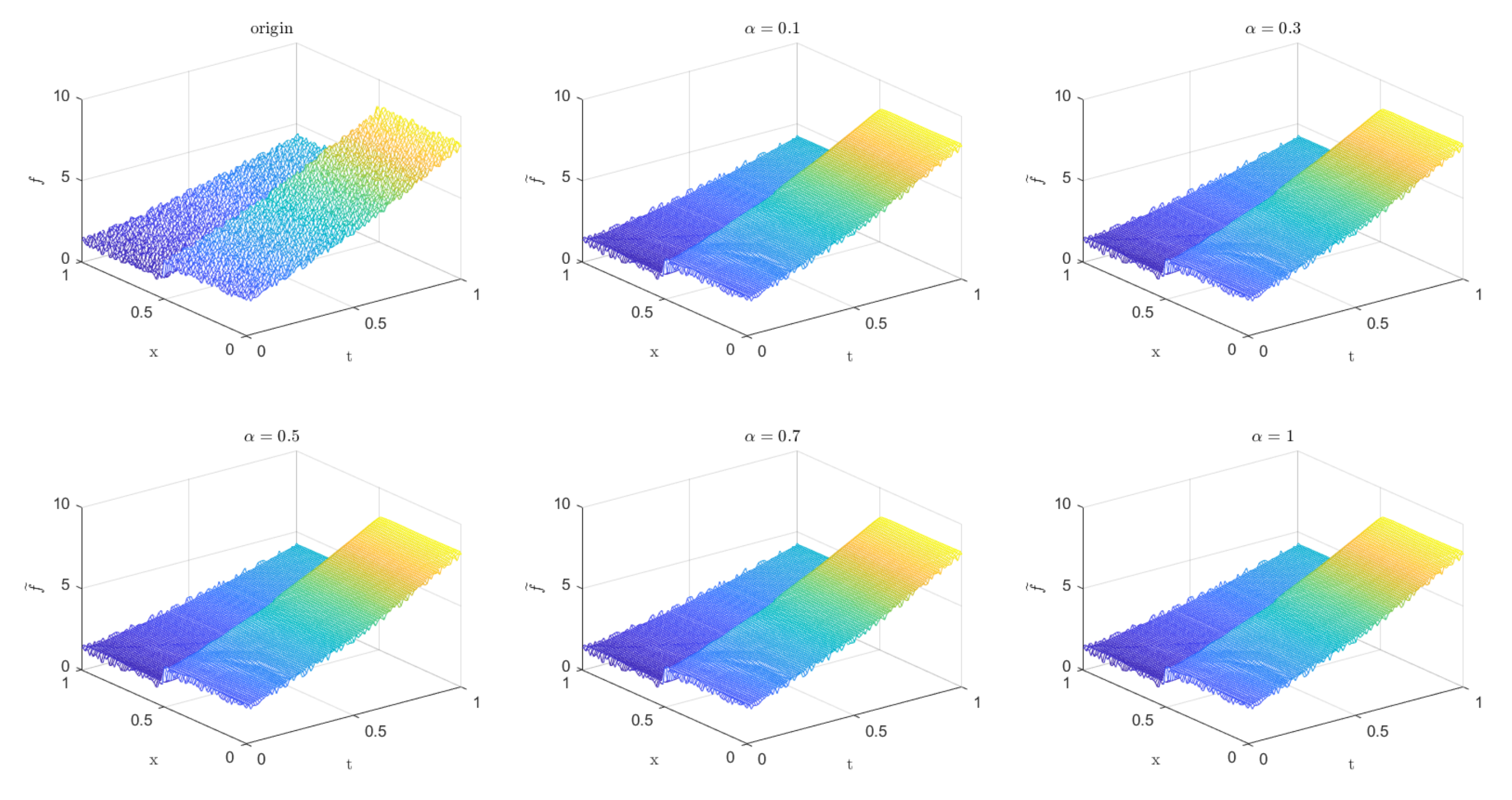}
	\caption{\textbf{Discover $\widetilde{f}$ for $h_{x}=h_{t}=1/100$ with $50\%$ uniform  noise  in Type 1  for Example \ref{ex2}.}\label{fig6.1.1.1}}
\end{figure*}

\begin{figure*}[!h]
	\centering
	\includegraphics[scale=0.35]{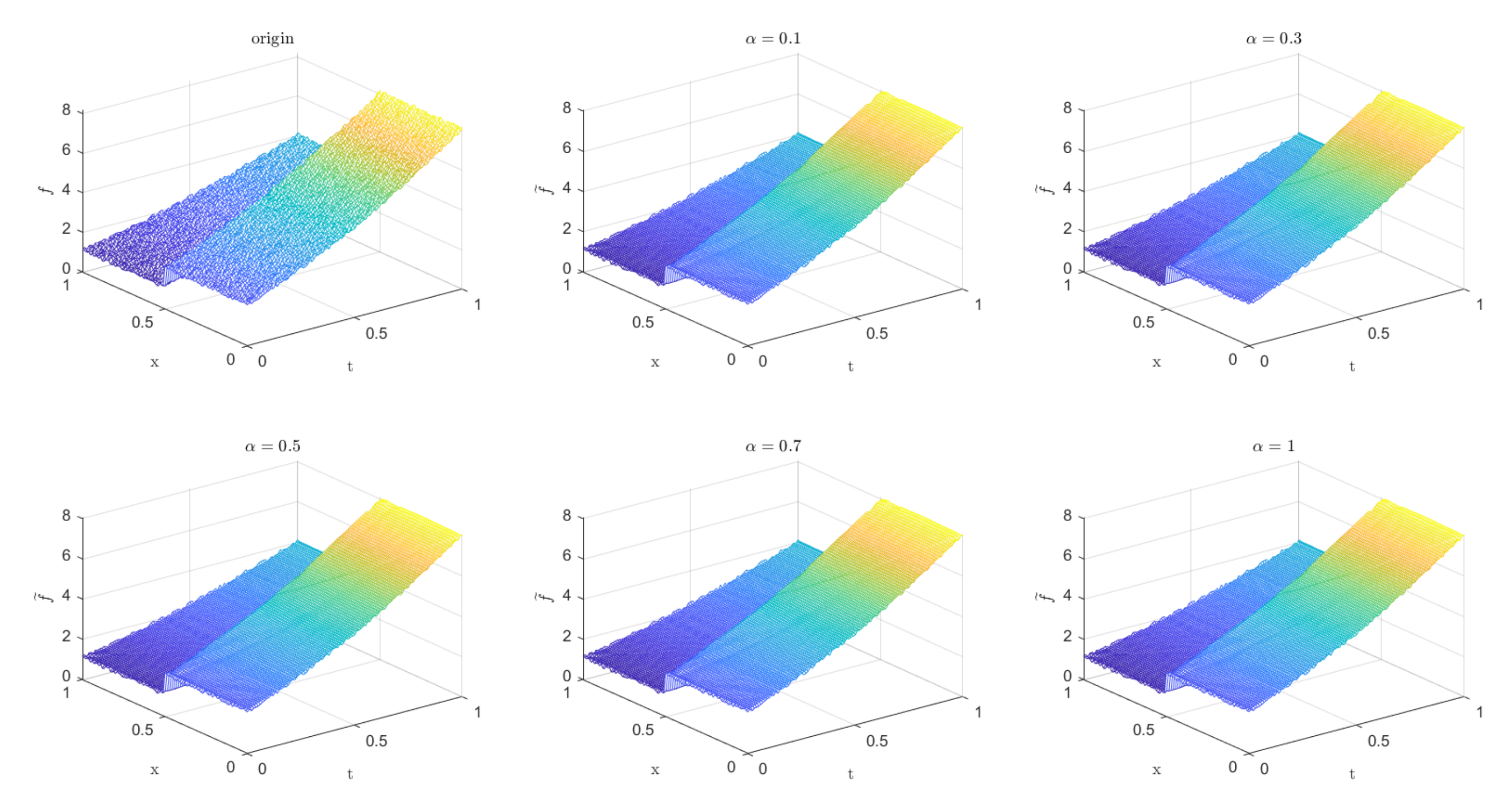}
	\caption{\textbf{Discover $\widetilde{f}$ for $h_{x}=h_{t}=1/100$ with $20\%$ uniform  noise  in Type 1  for Example \ref{ex2}.}\label{fig6.1.1.2}}
\end{figure*}

\begin{figure*}[!h]
	\centering
	\includegraphics[scale=0.35]{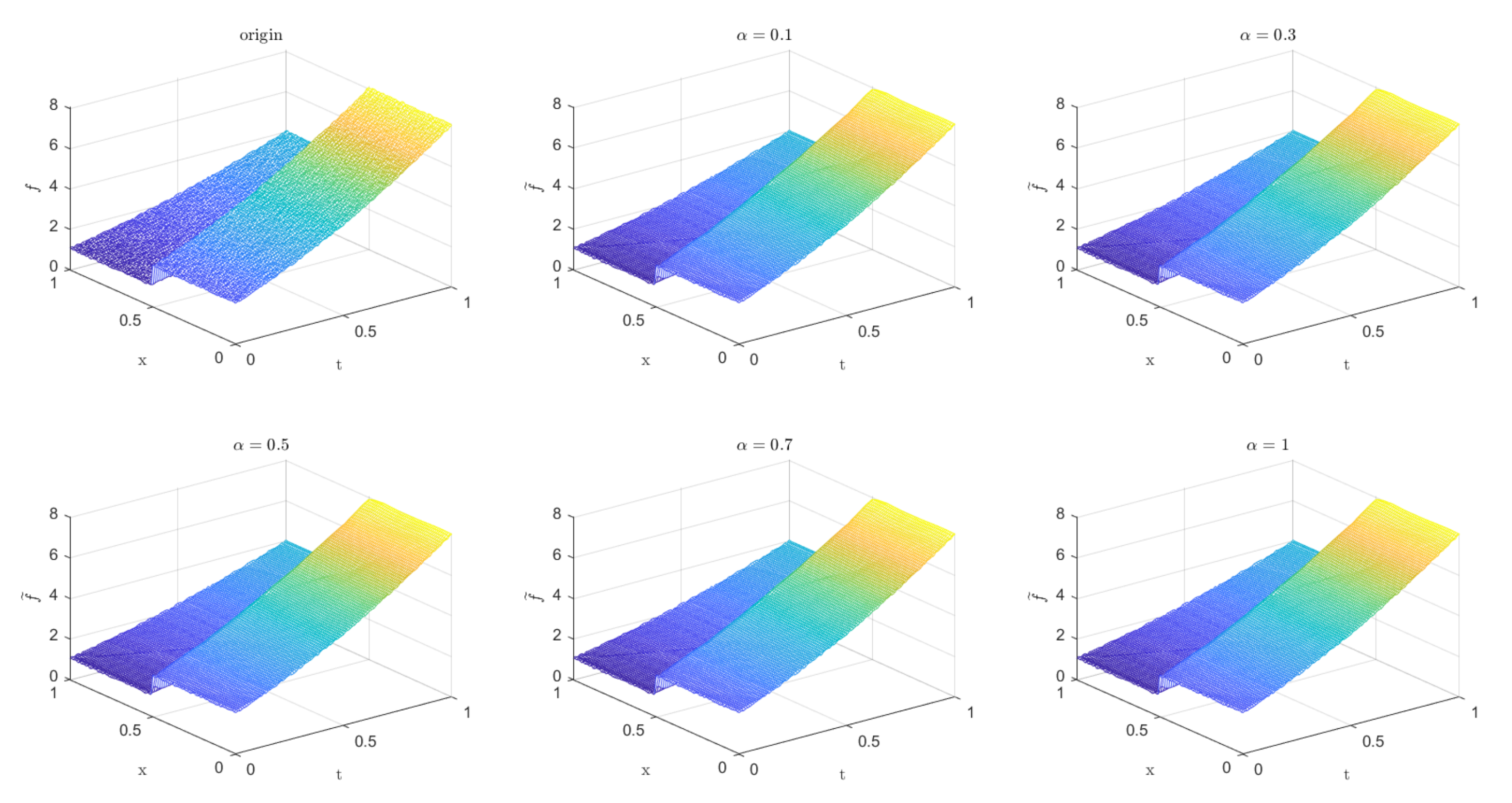}
	\caption{\textbf{Discover $\widetilde{f}$ for $h_{x}=h_{t}=1/100$ with $10\%$ uniform  noise  in Type 1  for Example \ref{ex2}.}\label{fig6.1.1.3}}
\end{figure*}

\begin{figure*}[!h]
	\centering
	\includegraphics[scale=0.35]{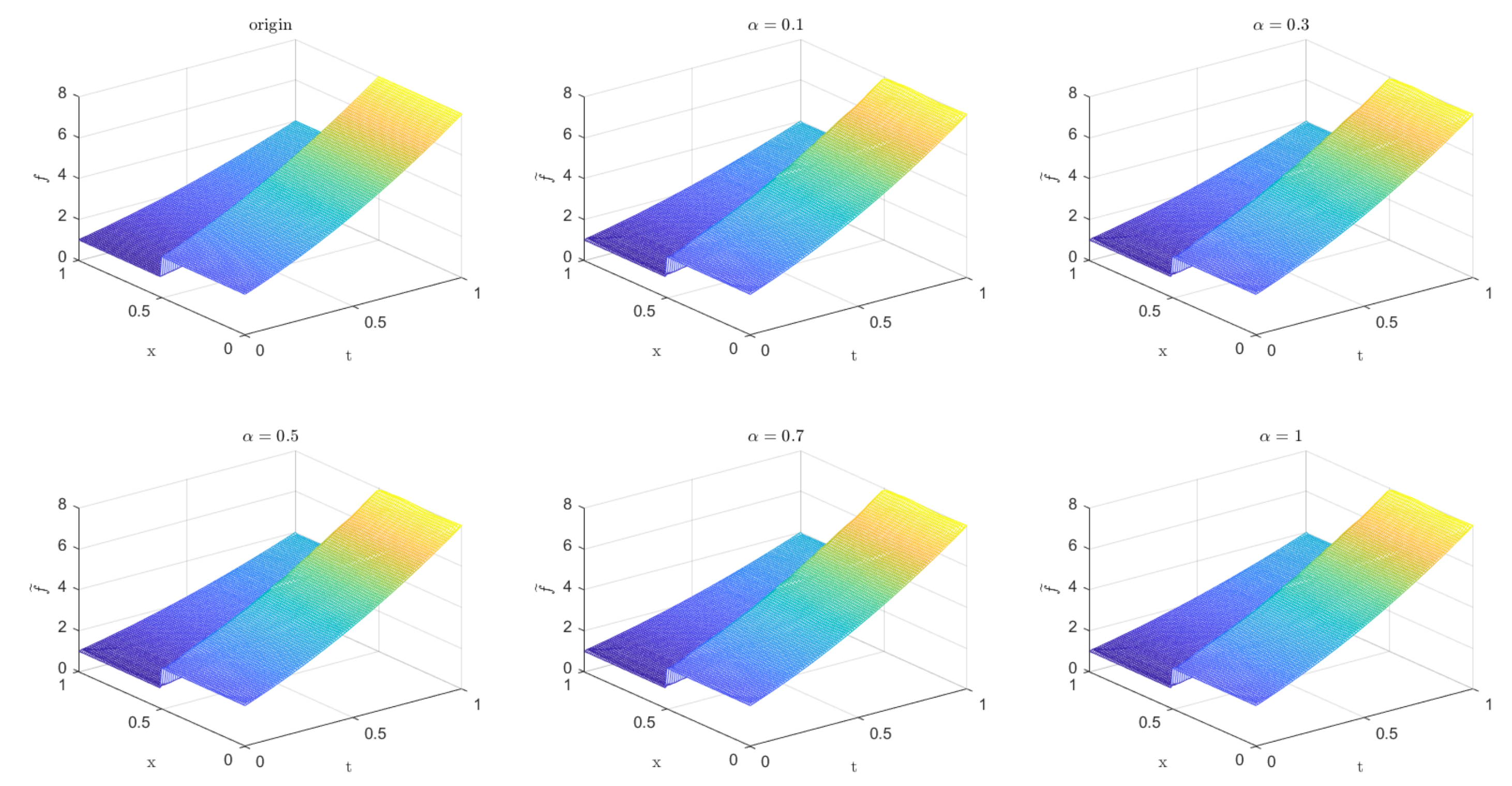}
	\caption{\textbf{Discover $\widetilde{f}$ for $h_{x}=h_{t}=1/100$ with clean data  in Type 1  for Example \ref{ex2}.}\label{fig6.1.1.4}}
\end{figure*}

Table \ref{Tab1} and Figures \ref{fig6.1.1.1}-\ref{fig6.1.1.4} show the relative error \eqref{D6.1} between the source function $f$ in Example \ref{ex2} and discovery $\widetilde{f}$ in  \eqref{C6} with epoch=$250$
 in Type 1.
From Table \ref{Tab1}, it can be found that our proposed algorithm is stable and accurate faced with different discretization sizes and noise level,
which is robust  even for  $50\%$  uniformly distributed  noise level.

\begin{table}[!ht]
	\centering
	\renewcommand\arraystretch{1.5}
	\caption{\textbf{The relative error $||e_{\hat{f}}||_{r}$  with uniform  noise in Type 1  for Example \ref{ex2}. }}\label{Tab1}
	\label{tab:the error comparison}
	\resizebox{\textwidth}{!}{
	\begin{tabular}{|c|c|ccccc|}
		\hline
	Threshold& \diagbox{$h_{x}$}{$\alpha$} & 0.1 & 0.3 & 0.5  &0.7 &1\\
    \hline
\multirow{3}{*}{$50\%$ noise}&1/25 &1.3898e-01
&1.3898e-01&1.3898e-01&1.3898e-01 &1.3898e-01\\
		\cline{2-7}
		&1/50 &8.1201e-02&8.1201e-02&8.1201e-02&8.1201e-02&	8.1201e-02\\
		\cline{2-7}
	   &1/100&3.4951e-02&3.4951e-02&3.4951e-02&3.4951e-02	&3.4951e-02\\
	   \hline
\multirow{3}{*}{$20\%$ noise}&1/25 &1.4190e-01&1.4190e-01	&1.4190e-01&1.4190e-01&1.4190e-01\\
		\cline{2-7}
		&1/50 &8.2515e-02&8.2515e-02&8.2515e-02&8.2515e-02	&8.2515e-02\\
		\cline{2-7}
	   &1/100&1.6557e-02&1.6557e-02&1.6557e-02&1.6557e-02	&1.6557e-02\\
	   \hline
\multirow{3}{*}{$10\%$ noise}&1/25 &1.3998e-01&1.3998e-01	&1.3998e-01&1.3998e-01&1.3998e-01\\
		\cline{2-7}
		&1/50 &8.1266e-02&8.1266e-02&8.1266e-02&8.1266e-02&	8.1266e-02\\
		\cline{2-7}
	   &1/100&8.1891e-03&8.1891e-03&8.1891e-03&8.1891e-03&8.1891e-03\\
	   \hline
\multirow{3}{*}{\tabincell{c}{clean data}}& 1/25 &1.2192e-01
&1.2192e-01&1.2192e-01&1.2192e-01 &1.2192e-01\\
		\cline{2-7}
		&1/50 &8.1341e-02
&8.1341e-02&8.1341e-02&8.1341e-02&8.1341e-02\\
		\cline{2-7}
	   &1/100&3.1629e-03&3.1629e-03&3.1629e-03&3.1629e-03&3.1629e-03\\
\hline
	\end{tabular}}
\end{table}

\subsubsection{Type 2 in DNN}
\begin{figure*}[!h]
	\centering
	\includegraphics[scale=0.35]{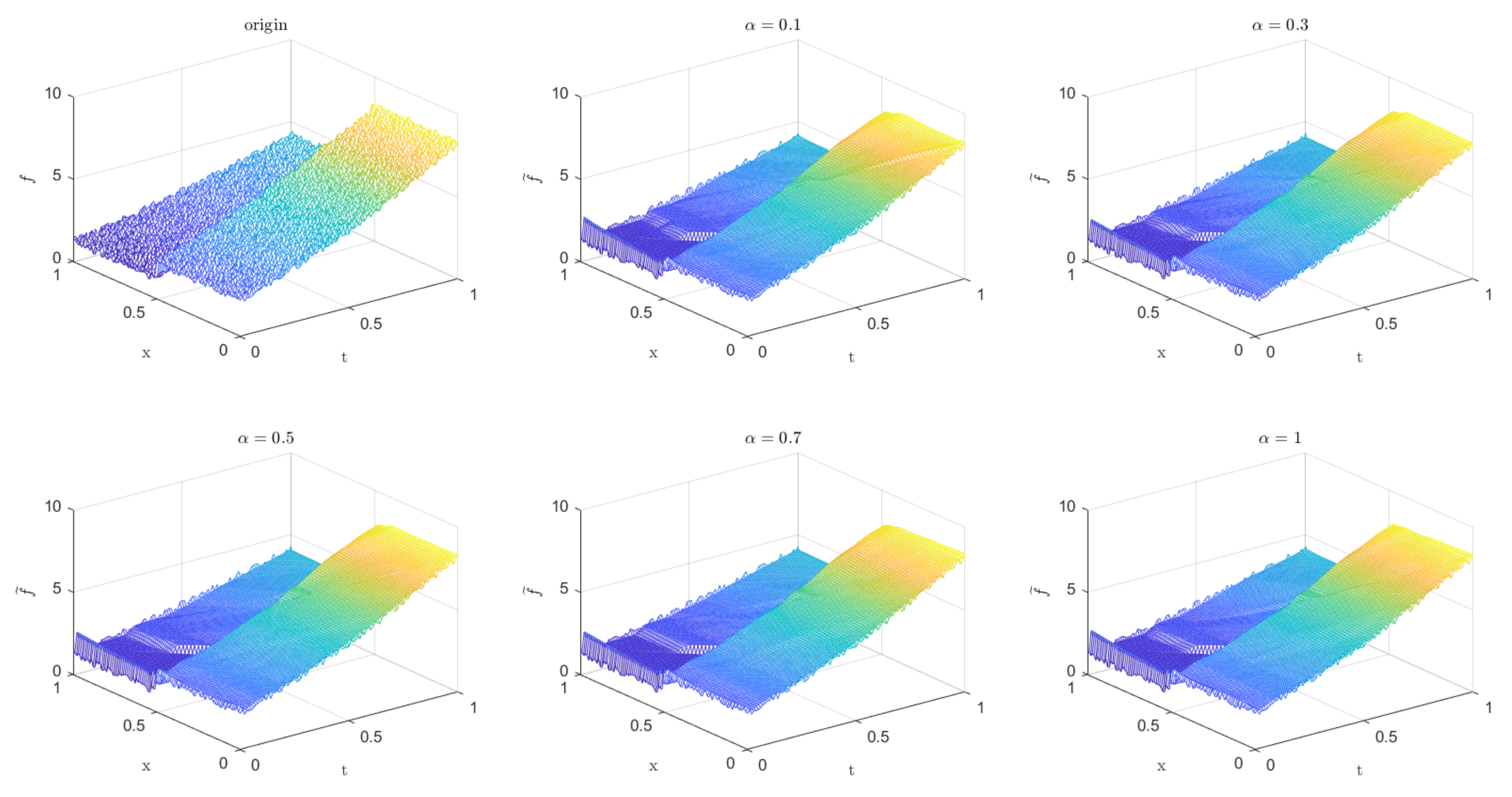}
	\caption{\textbf{Discover  $\widetilde{f}$ for $h_{x}=h_{t}=1/100$ with  $50\%$ uniform  noise in Type 2   for Example \ref{ex2}.}}\label{Fig6.1.2.1}
\end{figure*}

\begin{figure*}[!h]
	\centering
	\includegraphics[scale=0.35]{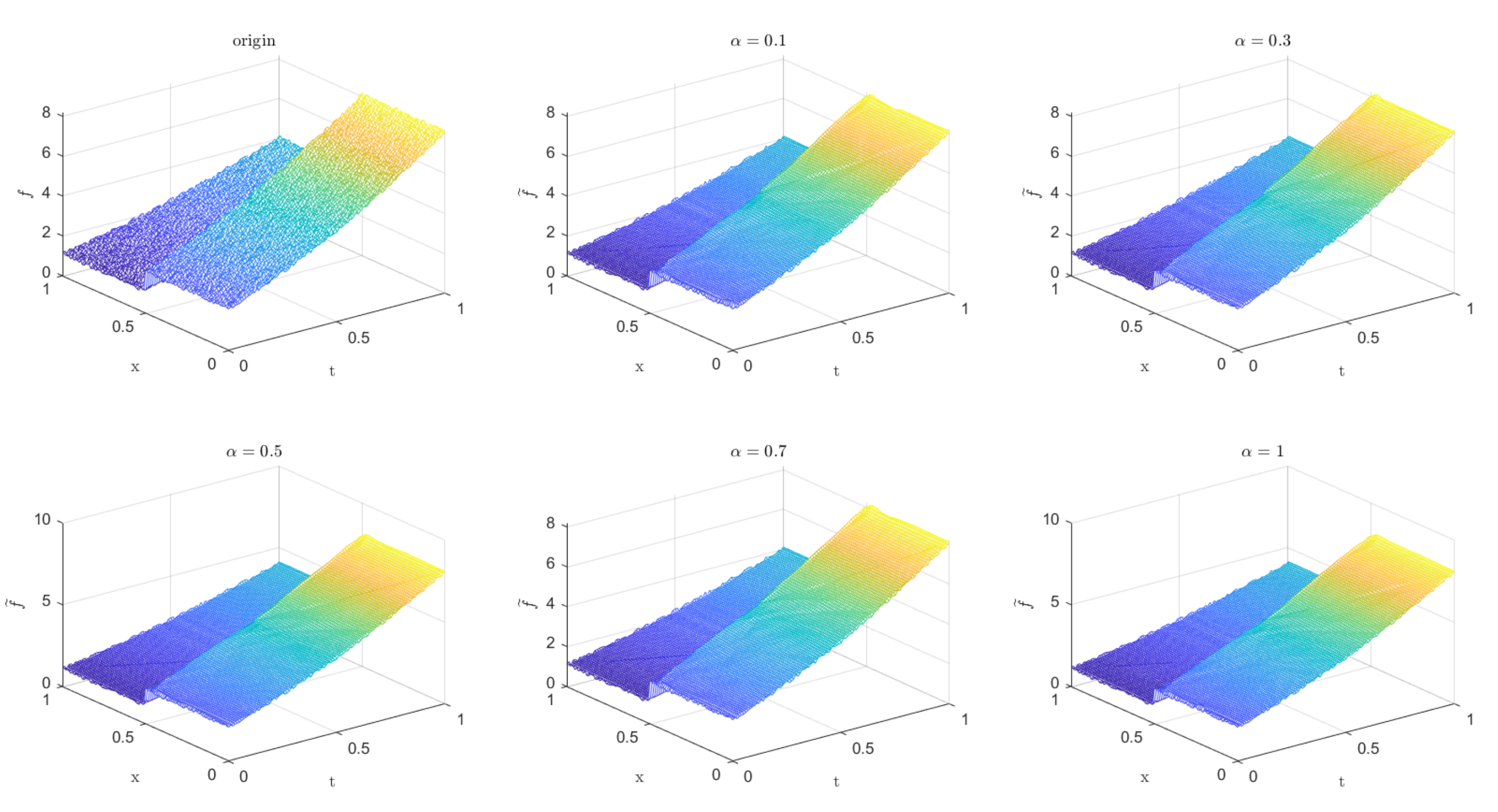}
	\caption{\textbf{Discover  $\widetilde{f}$ for $h_{x}=h_{t}=1/100$ with $20\%$ uniform  noise in Type 2   for Example \ref{ex2}.}}\label{Fig6.1.2.2}
\end{figure*}

\begin{figure*}[!h]
	\centering
	\includegraphics[scale=0.35]{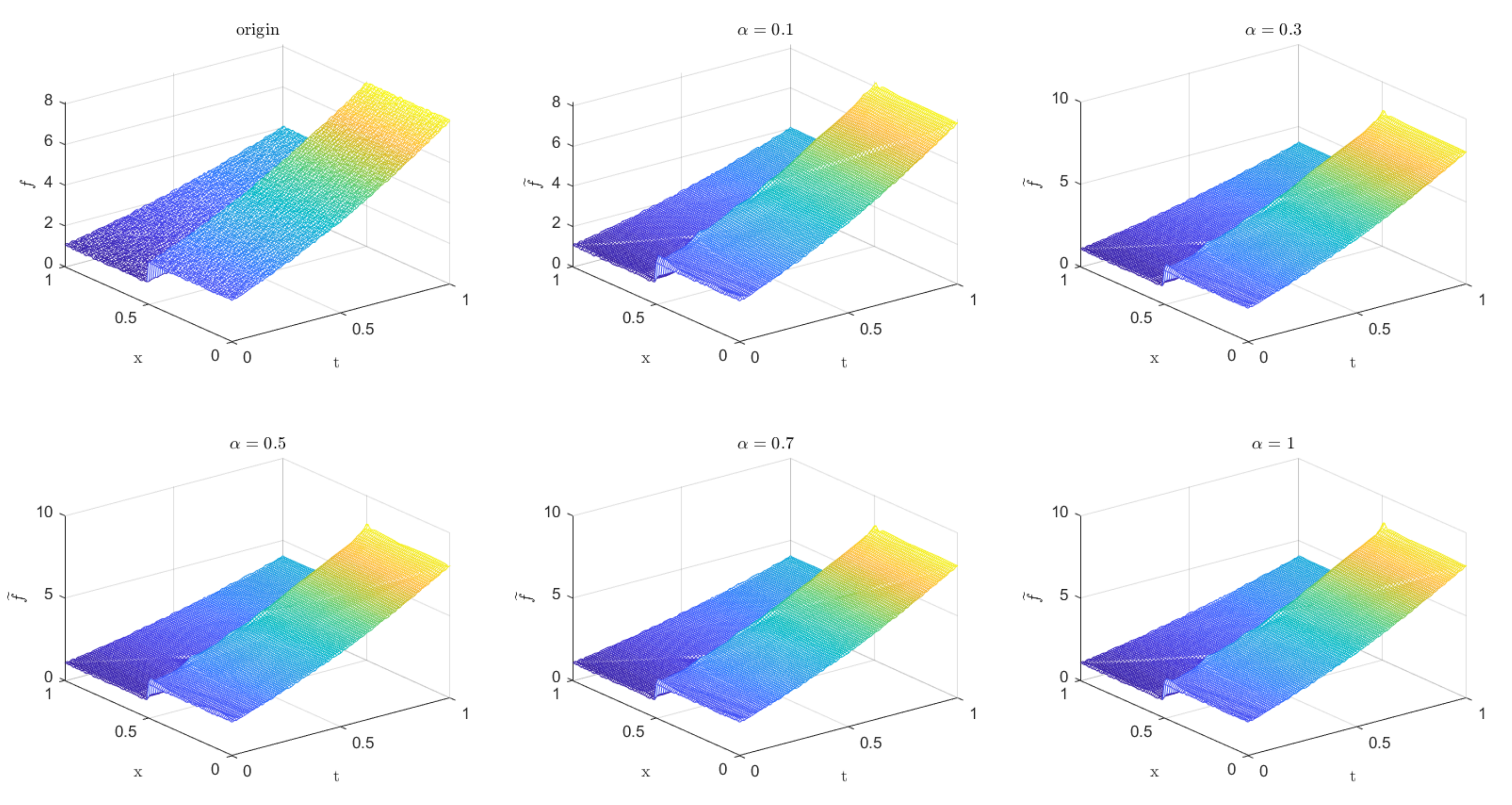}
	\caption{\textbf{Discover  $\widetilde{f}$ for $h_{x}=h_{t}=1/100$ with $10\%$ uniform  noise in Type 2   for Example \ref{ex2}.}}\label{Fig6.1.2.3}
\end{figure*}

\begin{figure*}[!h]
	\centering
	\includegraphics[scale=0.35]{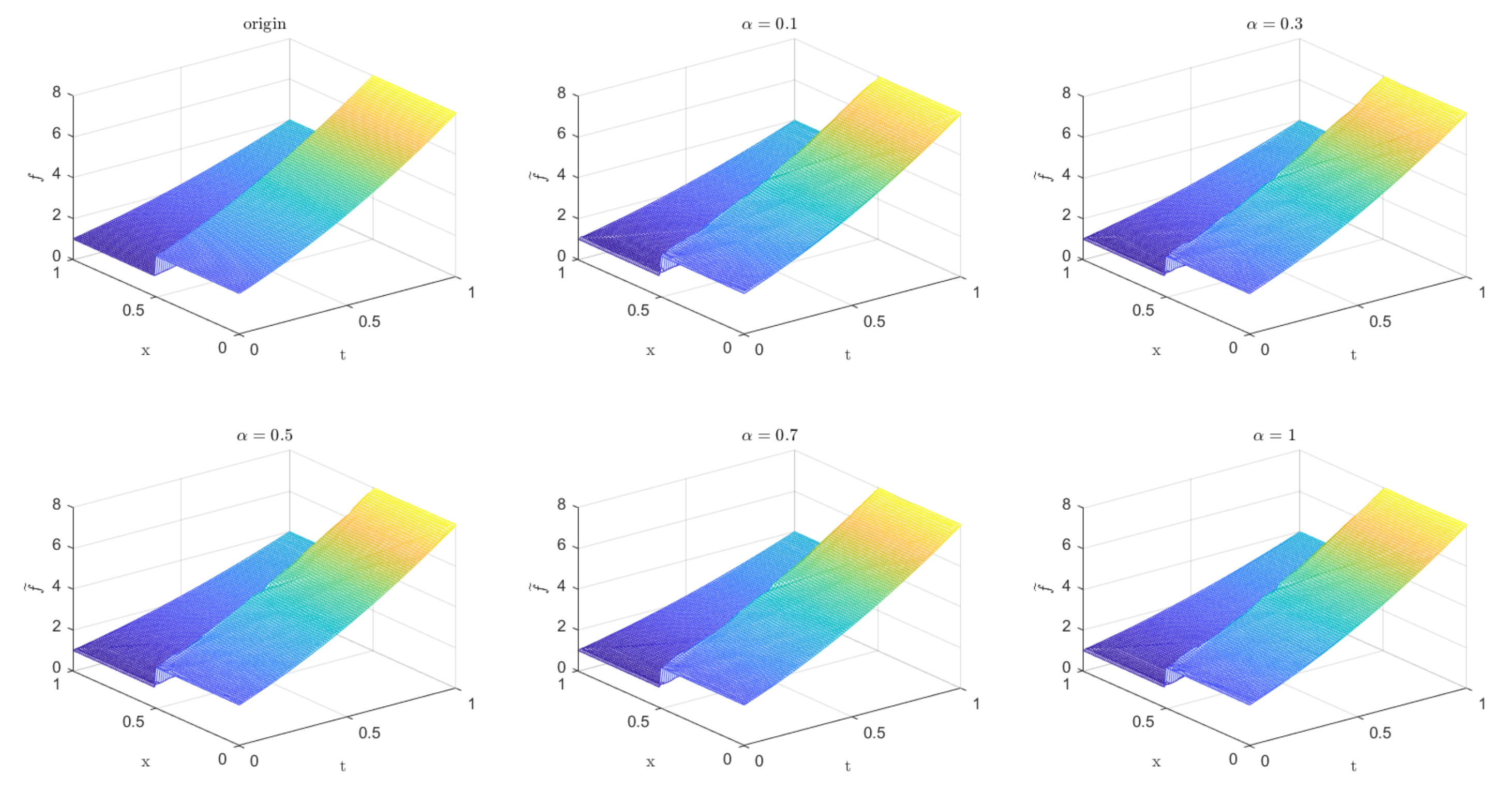}
	\caption{\textbf{Discover  $\widetilde{f}$ for $h_{x}=h_{t}=1/100$ with clean data in Type 2   for Example \ref{ex2}.}}\label{Fig6.1.2.4}
\end{figure*}

Table  \ref{Tab2} and Figures \ref{Fig6.1.2.1}-\ref{Fig6.1.2.4} show the relative error \eqref{D6.1} between the source function $f$ in Example \ref{ex2} and discovery $\widetilde{f}$ in  \eqref{C6}
with epoch=$250$  in Type 2.  The numerical experiments  are given to illustrate  the availability using deep learning.
In fact, from Table \ref{Tab2}, it can be found that our proposed algorithm is stable and accurate faced with different discretization sizes and noise level,
and is also robust  even for  $50\%$  uniformly distributed  noise level.

\begin{table}[!h]
	\centering
	\renewcommand\arraystretch{1.5}
	\caption{\textbf{The relative error  $||e_{\hat{f}}||_{r}$ with uniform  noise  in Type 2  for Example \ref{ex2}. }} \label{Tab2}
	\label{tab:the error comparison}
	\resizebox{\textwidth}{!}{
	\begin{tabular}{|c|c|ccccc|}
		\hline
	Threshold& \diagbox{$h_{x}$}{$\alpha$} & 0.1   &0.3  &0.5 &0.7 &1\\
	\hline
		    \multirow{3}{*}{$50\%$ noise}&1/25 &1.0774e-01&	1.0659e-01&1.0676e-01&1.0695e-01&1.0952e-01\\

		\cline{2-7}
		&1/50 &7.9438e-02&7.9069e-02&7.9516e-02&7.9440e-02&	8.1030e-02\\

		\cline{2-7}
	   &1/100&5.2110e-02&5.1856e-02&5.1991e-02&5.2098e-02	&5.2704e-02\\
	   \hline
	   	    \multirow{3}{*}{$20\%$ noise}&1/25 &1.1791e-01&	1.1925e-01&1.1919e-01&1.1859e-01&1.1716e-01\\

		\cline{2-7}
		&1/50 &7.9669e-02&8.0419e-02&8.0401e-02&7.9954e-02	&7.9223e-02\\

		\cline{2-7}
	   &1/100&1.5546e-02&1.5030e-02&1.5171e-02&1.4991e-02	&1.6229e-02\\
	   \hline
	   	   \multirow{3}{*}{$10\%$ noise} &1/25 &1.2006e-01&	1.2129e-01&1.2091e-01&1.2068e-01&1.2047e-01\\
		\cline{2-7}
		&1/50 &8.0694e-02&8.1230e-02&8.1086e-02&8.0716e-02&	8.0783e-02\\
		\cline{2-7}
	   &1/100&1.0260e-02&9.6543e-03&9.5003e-03&9.5757e-03&	1.0225e-02\\
	   \hline
	    \multirow{3}{*}{clean data}&1/25 &1.2154e-01&1.2198e-01&1.2206e-01	&1.2179e-01&1.2169e-01\\

		\cline{2-7}
		&1/50 &8.1152e-02&8.1369e-02&8.1410e-02&8.1255e-02&8.1141e-02\\

		\cline{2-7}
	   &1/100&4.8904e-03&3.7653e-03&3.7053e-03&3.6545e-03&4.7890e-03\\
	   \hline
	\end{tabular}}
\end{table}

\begin{example}\label{ex21}
Let us consider the following subdiffusion problem \eqref{ad1.1} with uniform  noise
\begin{equation*}
	\begin{cases}
		{_{0}^CD_t^{\alpha}u(x,t)}-\Delta u =f, \quad x \in [0,1], t \in [0,1];\\
		u(x,0)=\sqrt{x(1-x)};\\
		u(0,t)=0, u(1,t)=0,	
	\end{cases}
\end{equation*}
where the source function with the random noise $\delta$ is given by
\begin{equation*}
    f(x,t)=\begin{cases}
        (t+1)^{1/4}\left(1+\chi_{(0,1/2)}(x)\right)+\delta,\quad 0\le x \le 1/2;\\
        (t+1)^2\left(1+\chi_{(0,1/2)}(x)\right)+\delta,\quad {\rm otherwise}.
    \end{cases}
\end{equation*}
Here $\delta$ is the  uniform  noise with $10\%$ (clean data) $10\%$, $20\%$, $50\%$ level for a single uniform distributed random number in the interval $(0,1)$, respectively.
In this experiment, we train  the deep learning discovery $\widetilde{f}$ in \eqref{C6}.
\end{example}
\subsubsection{Type 1 in DNN}
\begin{figure*}[!h]
	\centering

	\includegraphics[scale=0.35]{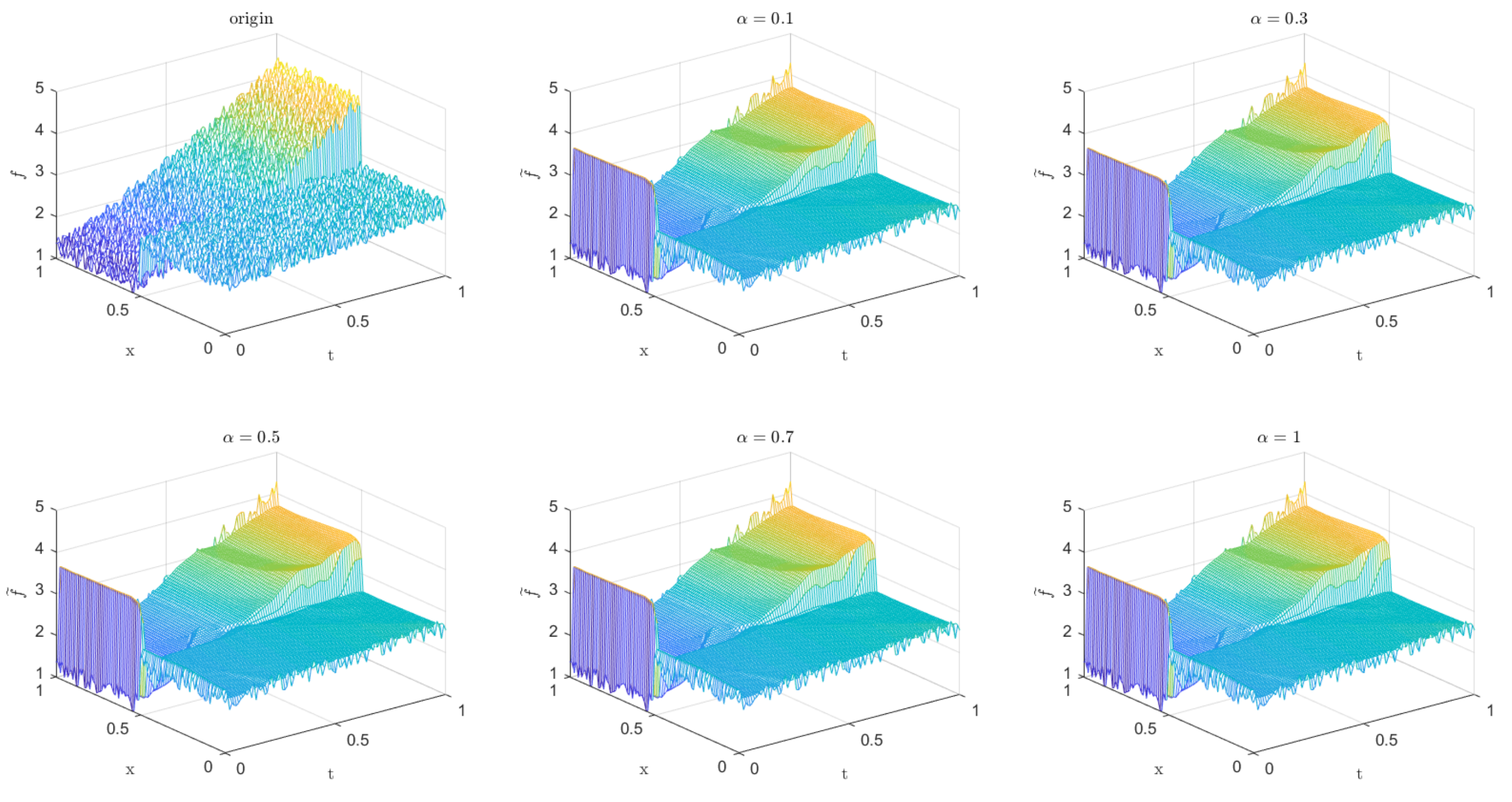}
	\caption{\textbf{Discover $\widetilde{f}$ for $h_{x}=h_{t}=1/100$ with $50\%$ uniform  noise in Type 1 for Example \ref{ex21}.}\label{T3.2.1.1}}
\end{figure*}

\begin{figure*}[!h]
	\centering

	\includegraphics[scale=0.35]{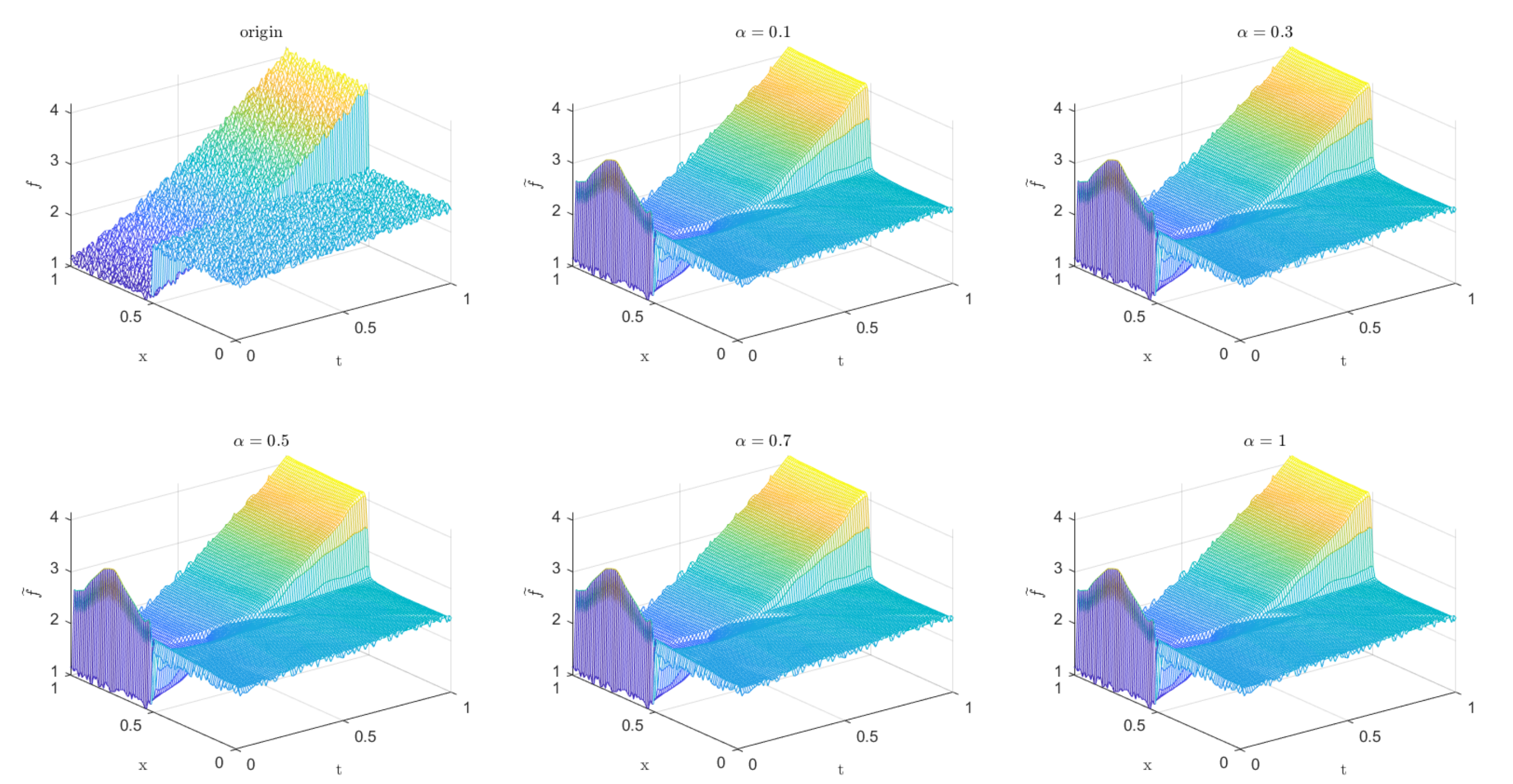}
	\caption{\textbf{Discover $\widetilde{f}$ for $h_{x}=h_{t}=1/100$ with $20\%$ uniform  noise in Type 1 for Example \ref{ex21}.}\label{T3.2.1.2}}
\end{figure*}

\begin{figure*}[!h]
	\centering

	\includegraphics[scale=0.35]{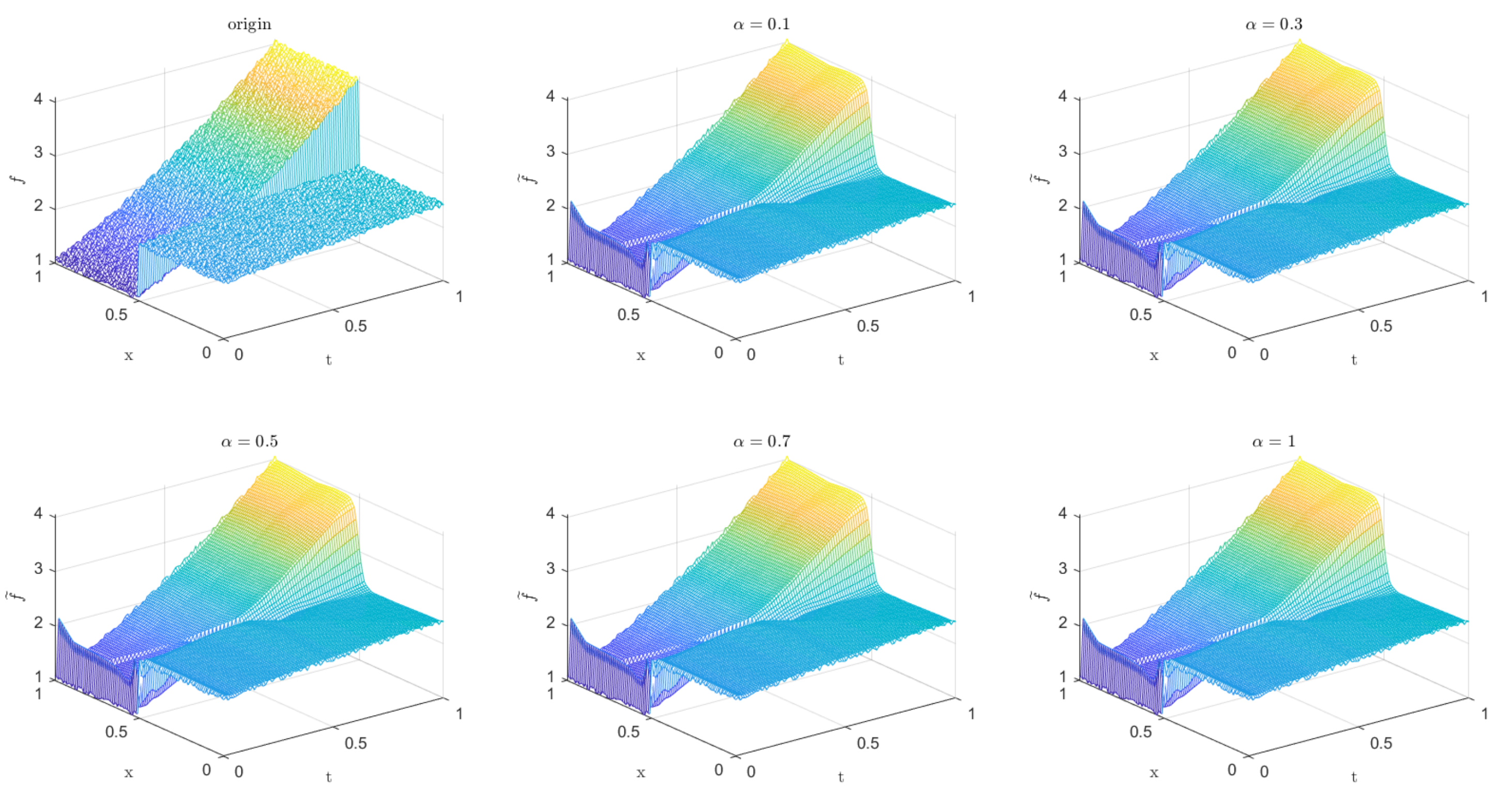}
	\caption{\textbf{Discover $\widetilde{f}$ for $h_{x}=h_{t}=1/100$ with $10\%$ uniform  noise in Type 1 for Example \ref{ex21}.}\label{T3.2.1.3}}
\end{figure*}

\begin{figure*}[!h]
	\centering

	\includegraphics[scale=0.35]{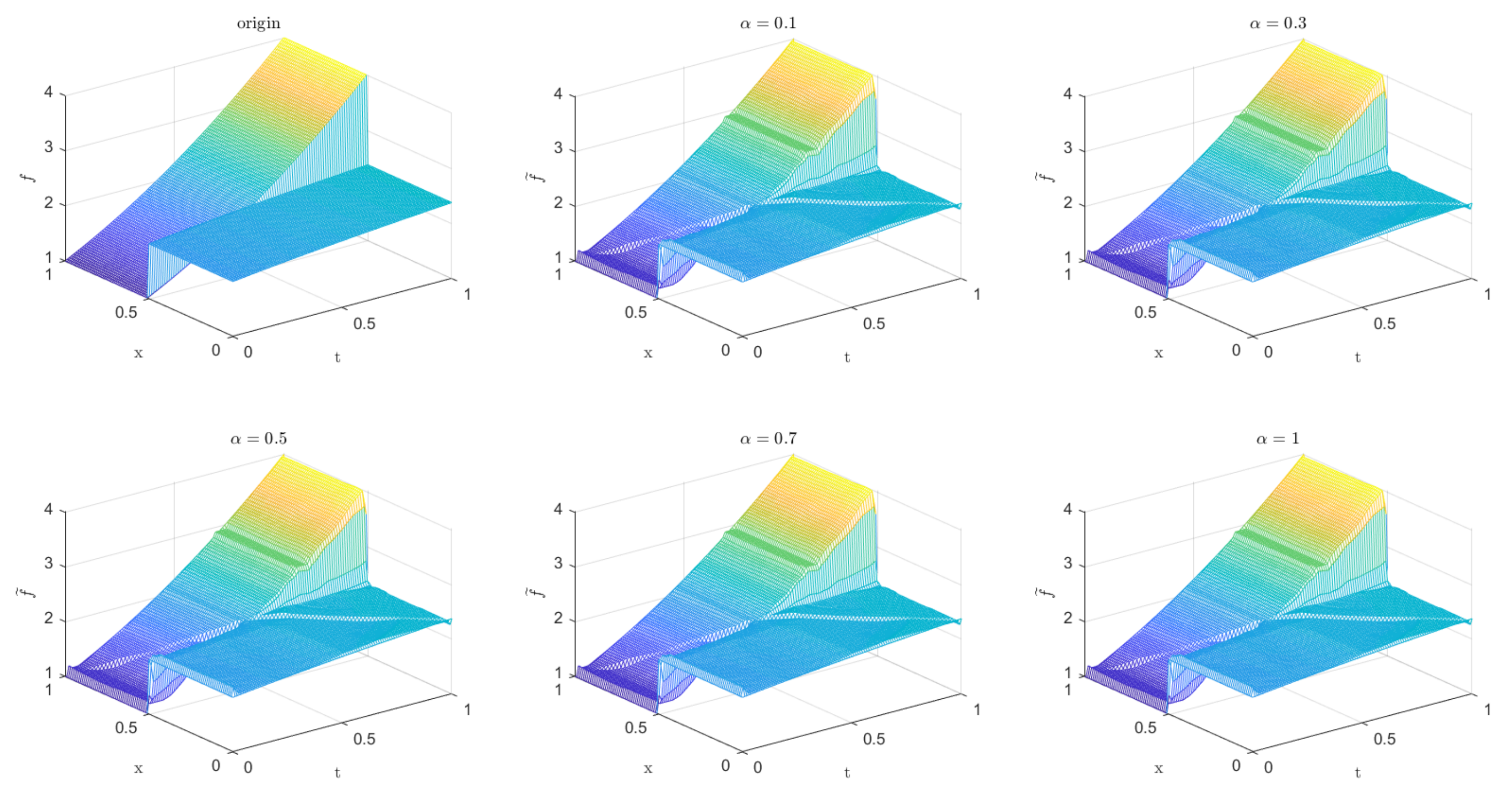}
	\caption{\textbf{Discover $\widetilde{f}$ for $h_{x}=h_{t}=1/100$ with clean data in Type 1 for Example \ref{ex21}.}\label{T3.2.1.4}}
\end{figure*}

Table  \ref{Tab3} and Figures \ref{T3.2.1.1}-\ref{T3.2.1.4} show the relative error \eqref{D6.1} between the source function $f$ in Example \ref{ex21} and discovery $\widetilde{f}$ in  \eqref{C6}
with different noise levels and epoch=$200$ in Type 1. The numerical experiments  are given to illustrate  the availability using deep learning for the different noise levels.
From Table \ref{Tab3}, it can be found that our proposed algorithm is stable and accurate faced with different discretization sizes and noise level,
which is robust  even for  $50\%$  uniformly distributed  noise level.
From  Figures \ref{T3.2.1.1}-\ref{T3.2.1.4}, it appears a layer or blows up at $t=0$ (boundary noise pollution), since the low time regularity \cite{ShCh:2020,SOG:17}.

\begin{table}[!ht]
	\centering
	\renewcommand\arraystretch{1.5}
	\caption{\textbf{The relative error $||e_{\hat{f}}||_{r}$ with uniform  noise in Type 1 for Example \ref{ex21}.  }}\label{Tab3}
	\label{tab:the error comparison}
	\resizebox{\textwidth}{!}{
	\begin{tabular}{|c|c|ccccc|}
		\hline
	Noise Level& \diagbox{$h_{x}$}{$\alpha$} & 0.1 & 0.3 & 0.5  &0.7 &1\\
		   \hline
	    \multirow{3}{*}{\tabincell{c}{ $50\%$ noise \\}}&1/25 &1.5132e-01&1.5132e-01&1.5132e-01&1.5132e-01&	1.5132e-01\\
		\cline{2-7}
		&1/50 &1.2040e-01&1.2040e-01&1.2040e-01&1.2040e-01&1.2040e-01\\
		\cline{2-7}
	   &1/100&1.1607e-01&1.1607e-01&1.1607e-01&1.1607e-01&1.1607e-01\\
	   \hline
  	    \multirow{3}{*}{\tabincell{c}{$20\%$ noise \\}}&1/25 &1.1360e-01&1.1360e-01&1.1360e-01&1.1360e-01&	1.1360e-01\\
		\cline{2-7}
		&1/50 &8.3587e-02&8.3587e-02&8.3587e-02&8.3587e-02&8.3587e-02\\
		\cline{2-7}
	   &1/100&6.1966e-02&6.1966e-02&6.1966e-02&6.1966e-02&6.1966e-02\\
\hline
	    \multirow{3}{*}{\tabincell{c}{ $10\%$ noise\\}}&1/25 &5.8198e-02&5.8198e-02&5.8198e-02&5.8198e-02	&5.8198e-02\\
		\cline{2-7}
		&1/50 &4.2641e-02&4.2641e-02&4.2641e-02&4.2641e-02&4.2641e-02\\
		\cline{2-7}
	   &1/100&3.5011e-02&3.5011e-02&3.5011e-02&3.5011e-02&3.5011e-02	\\
	   \hline
	    \multirow{3}{*}{\tabincell{c}{ clean data\\}}&1/25 &5.2867e-02&5.2867e-02&5.2867e-02&5.2867e-02	&5.2867e-02\\
		\cline{2-7}
		&1/50 &3.7670e-02&3.7670e-02&3.7670e-02&3.7670e-02&3.7670e-02\\
		\cline{2-7}
	   &1/100&2.0318e-02&2.0318e-02&2.0318e-02&2.0318e-02&2.0318e-02	\\
	   \hline
	\end{tabular}}
\end{table}


\subsubsection{Type 2 in DNN}
\begin{figure*}[!ht]
	\centering
	\includegraphics[scale=0.35]{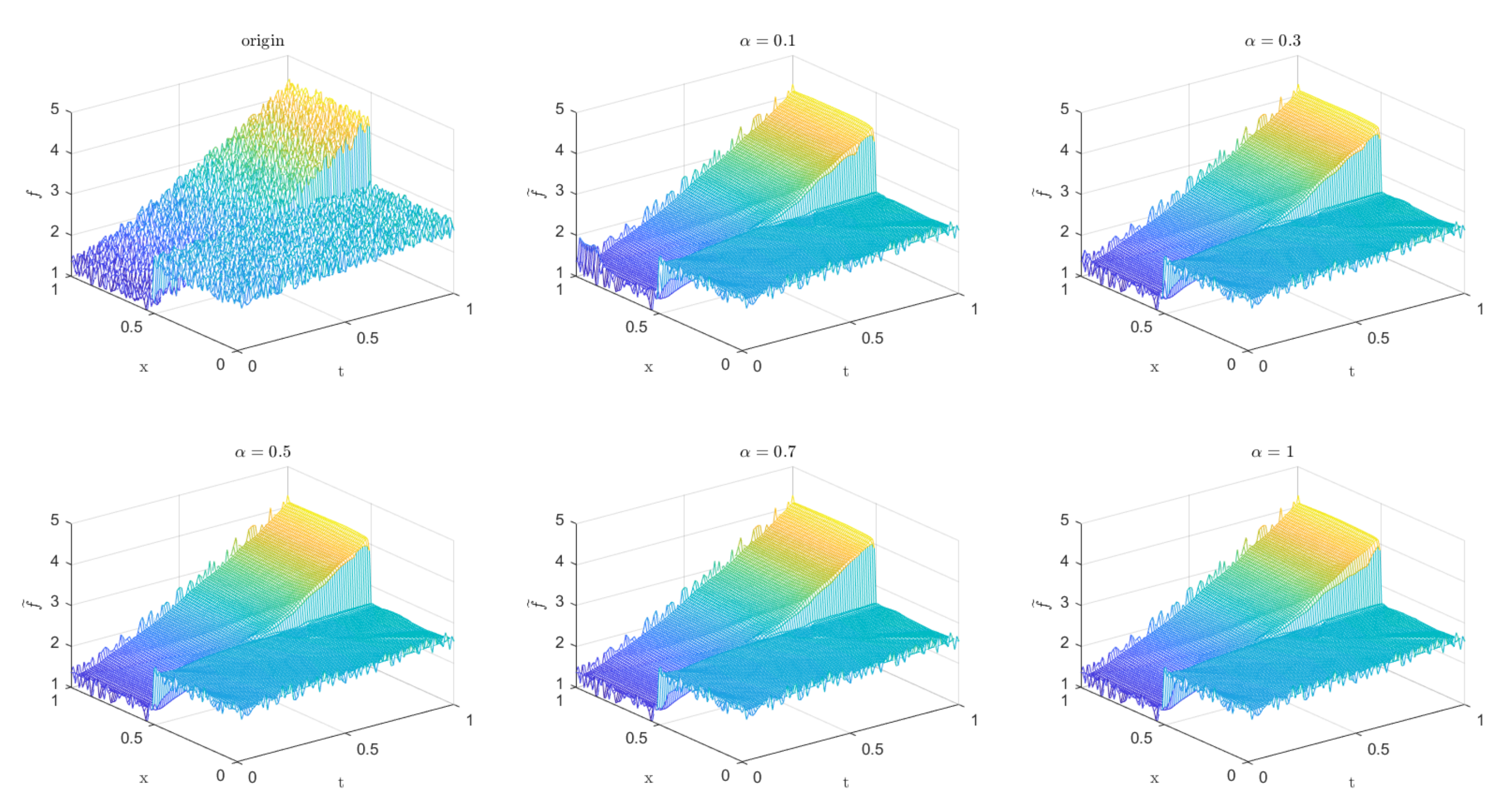}
	\caption{\textbf{Discover $\widetilde{f}$ for $h_{x}=h_{t}=1/100$ with $50\%$ uniform  noise in Type 2 for Example \ref{ex21}.  }}\label{T4.6.2.2.1}
\end{figure*}

\begin{figure*}[!h]
	\centering
	\includegraphics[scale=0.35]{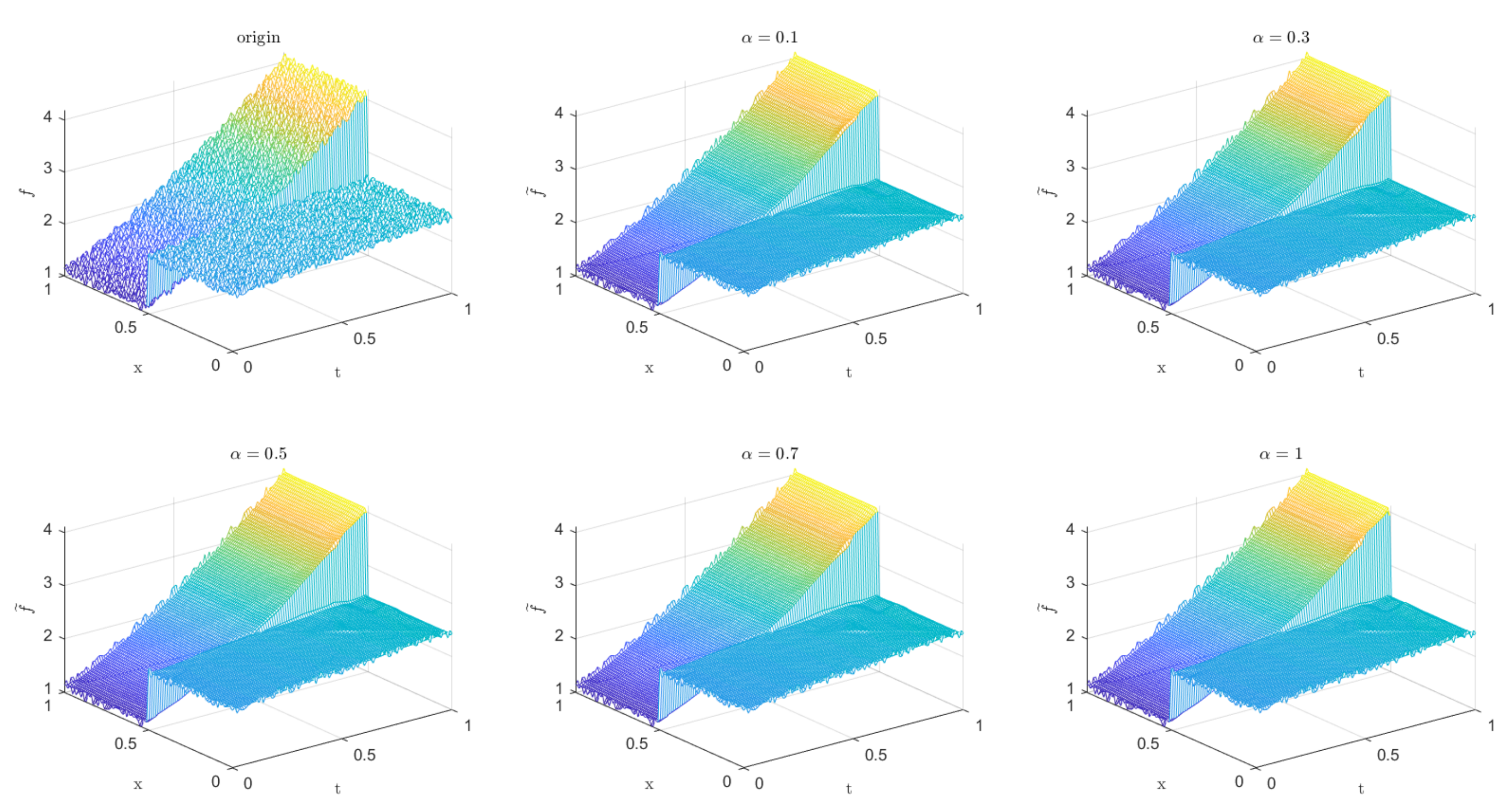}
	\caption{\textbf{Discover $\widetilde{f}$ for $h_{x}=h_{t}=1/100$ with $20\%$ uniform  noise in Type 2 for Example \ref{ex21}.  }}\label{T4.6.2.2.2}
\end{figure*}

\begin{figure*}[!h]
	\centering
	\includegraphics[scale=0.35]{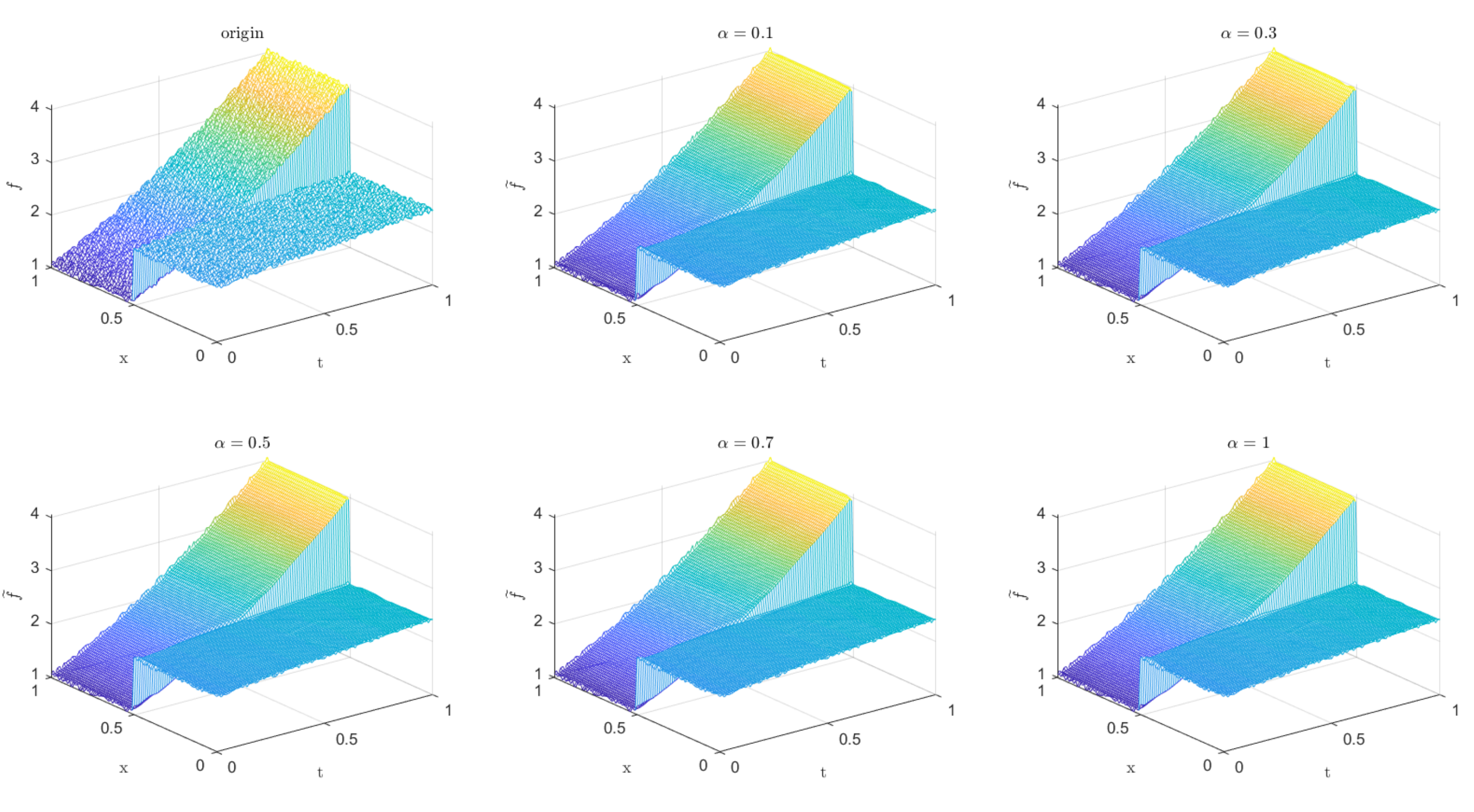}
	\caption{\textbf{Discover $\widetilde{f}$ for $h_{x}=h_{t}=1/100$ with $10\%$ uniform  noise in Type 2 for Example \ref{ex21}.  }}\label{T4.6.2.2.3}
\end{figure*}

\begin{figure*}[!h]
	\centering
	\includegraphics[scale=0.35]{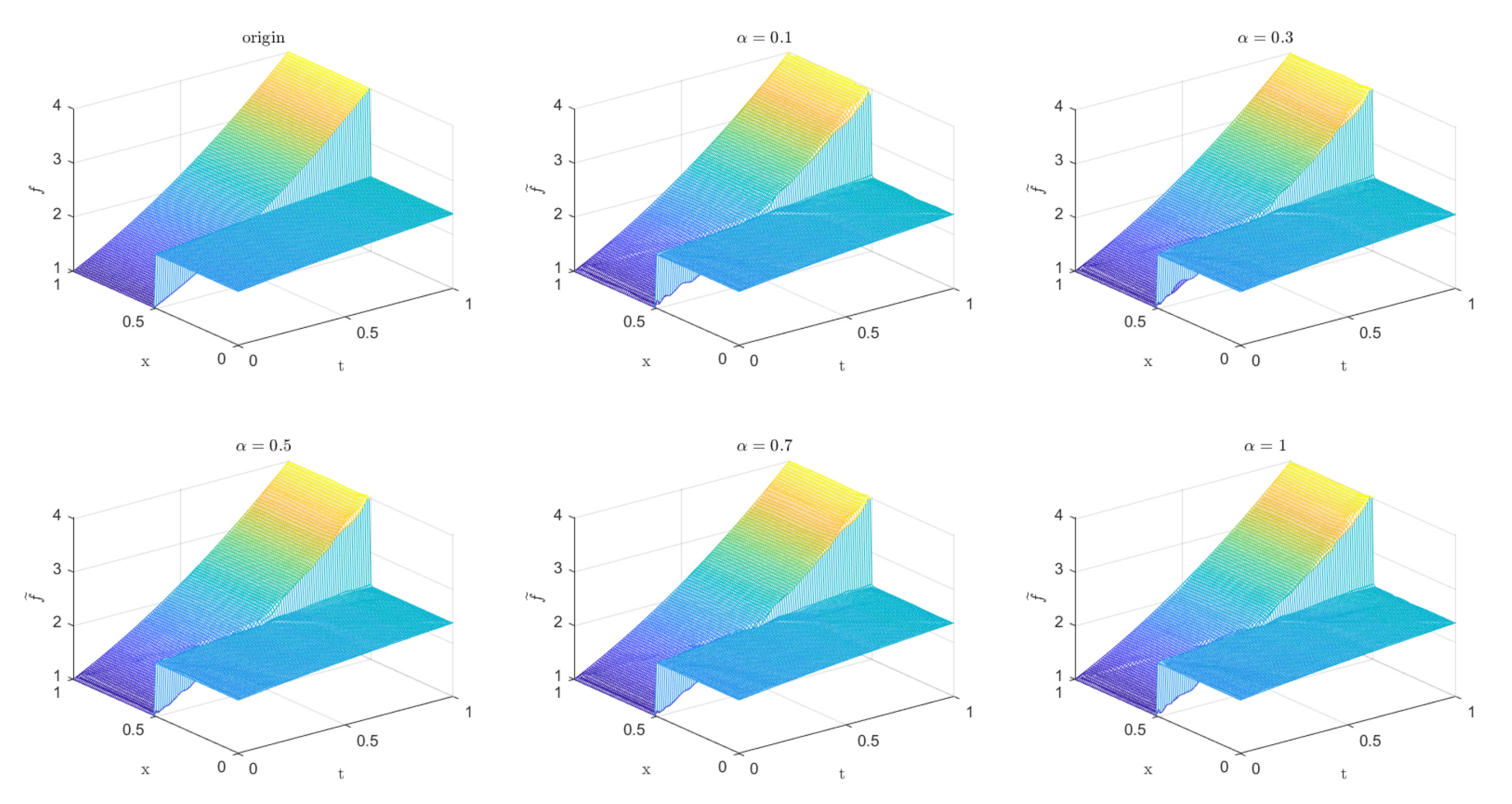}
	\caption{\textbf{Discover $\widetilde{f}$ for $h_{x}=h_{t}=1/100$ with clean data in Type 2 for Example \ref{ex21}.  }}\label{T4.6.2.2.4}
\end{figure*}

Table  \ref{Tab4} and  Figures  \ref{T4.6.2.2.1}-\ref{T4.6.2.2.4} show the relative error \eqref{D6.1} between the source function $f$ in Example \ref{ex2} and discovery $\widetilde{f}$ in  \eqref{C6}
with different noise levels and epoch=$270$ in Type 2. The numerical experiments  are given to illustrate  the availability using deep learning for the different noise levels.
From Table \ref{Tab4}, it can be found that our proposed algorithm is stable and accurate faced with different discretization sizes and noise level,
which is robust  even for  $50\%$  uniformly distributed  noise level.
From  Figures  \ref{T4.6.2.2.1}-\ref{T4.6.2.2.4}, it also appears a little layer or blows up at $t=0$ (boundary noise pollution), since the low time regularity \cite{ShCh:2020,SOG:17}.
It seems that Type 2 is better than Type1  from observed data.

\begin{table}[!ht]
	\centering
	\renewcommand\arraystretch{1.5}
	\caption{\textbf{The relative error  $||e_{\hat{f}}||_{r}$ with uniform  noise  in Type 2 for Example \ref{ex21}.   }} \label{Tab4}
	\label{tab:the error comparison}
	\resizebox{\textwidth}{!}{
	\begin{tabular}{|c|c|ccccc|}
		\hline
	Noise Level& \diagbox{$h_{x}$}{$\alpha$} & 0.1   &0.3  &0.5 &0.7 &1\\
	\hline
		   	 \multirow{3}{*}{\tabincell{c}{$50\%$ noise}}&1/25 &8.6685e-02&	8.5475e-02&8.4797e-02&8.4731e-02&8.1521e-02
\\
		\cline{2-7}
		&1/50 &7.2426e-02&7.1836e-02&7.1454e-02&7.1377e-02&	6.9543e-02\\
		\cline{2-7}
	   &1/100&6.2809e-02&6.2084e-02&6.1725e-02&6.1470e-02&	6.1917e-02\\
	   \hline
	   	   	 \multirow{3}{*}{\tabincell{c}{$20\%$ noise}}&1/25 &6.6977e-02&	6.7135e-02&6.6908e-02&6.7067e-02&6.7562e-02\\
		\cline{2-7}
		&1/50 &4.8547e-02&4.8697e-02&4.8523e-02&4.8578e-02&	4.8868e-02\\
		\cline{2-7}
	   &1/100&2.7690e-02&2.7786e-02&2.7643e-02&2.7814e-02&	2.7819e-02\\
	   \hline
	    \multirow{3}{*}{\tabincell{c}{$10\%$ noise}}&1/25 &6.3849e-02&6.3655e-02&	6.3809e-02&6.3935e-02&6.3953e-02\\
		\cline{2-7}
		&1/50 &4.3969e-02&4.3909e-02&4.3996e-02&4.4066e-02	&4.4010e-02\\
		\cline{2-7}
	   &1/100&1.8282e-02&1.8648e-02&1.8545e-02&1.8526e-02&	1.8353e-02\\
	   \hline
	    \multirow{3}{*}{\tabincell{c}{clean data}}&1/25 &5.6927e-02&5.8353e-02	&5.8272e-02&5.7993e-02&5.8502e-02\\
		\cline{2-7}
		&1/50 &3.8296e-02&3.9283e-02&3.9420e-02&3.9318e-02&	3.9337e-02\\
		\cline{2-7}
	   &1/100&5.4745e-03&4.4494e-03&4.9254e-03&4.7000e-03	&3.9269e-03\\
	   \hline
	\end{tabular}}
\end{table}

\subsection{Numerical Experiments for discovery with Guassian noise}
\begin{example}\label{ex31}
Let us consider the Example \ref{ex21} again with Guassian noise, namely,
\begin{equation*}
	\begin{cases}
		{_{0}^CD_t^{\alpha}u(x,t)}-\Delta u =f, \quad x \in [0,1], t \in [0,1];\\
		u(x,0)=\sqrt{x(1-x)};\\
		u(0,t)=0, u(1,t)=0,	
	\end{cases}
\end{equation*}
where the source function with the random noise $\eta$ is given by
\begin{equation*}
    f(x,t)=\begin{cases}
        (t+1)^{1/4}\left(1+\chi_{(0,1/2)}(x)\right)+\eta,\quad 0\le x \le 1/2;\\
        (t+1)^2\left(1+\chi_{(0,1/2)}(x)\right)+\eta,\quad {\rm otherwise}.
    \end{cases}
\end{equation*}
Here $\eta$ is the  Guassian  noise with $10\%$ (clean data) $10\%$, $20\%$, $50\%$ level, respectively, which is statistical noise having a probability density function equal
to that of the normal distribution.
In this experiment, we train  the deep learning discovery $\widetilde{f}$ in \eqref{C6}.
\end{example}

\subsubsection{Type 1 in DNN}
\begin{figure*}[!ht]
	\centering
	\includegraphics[scale=0.35]{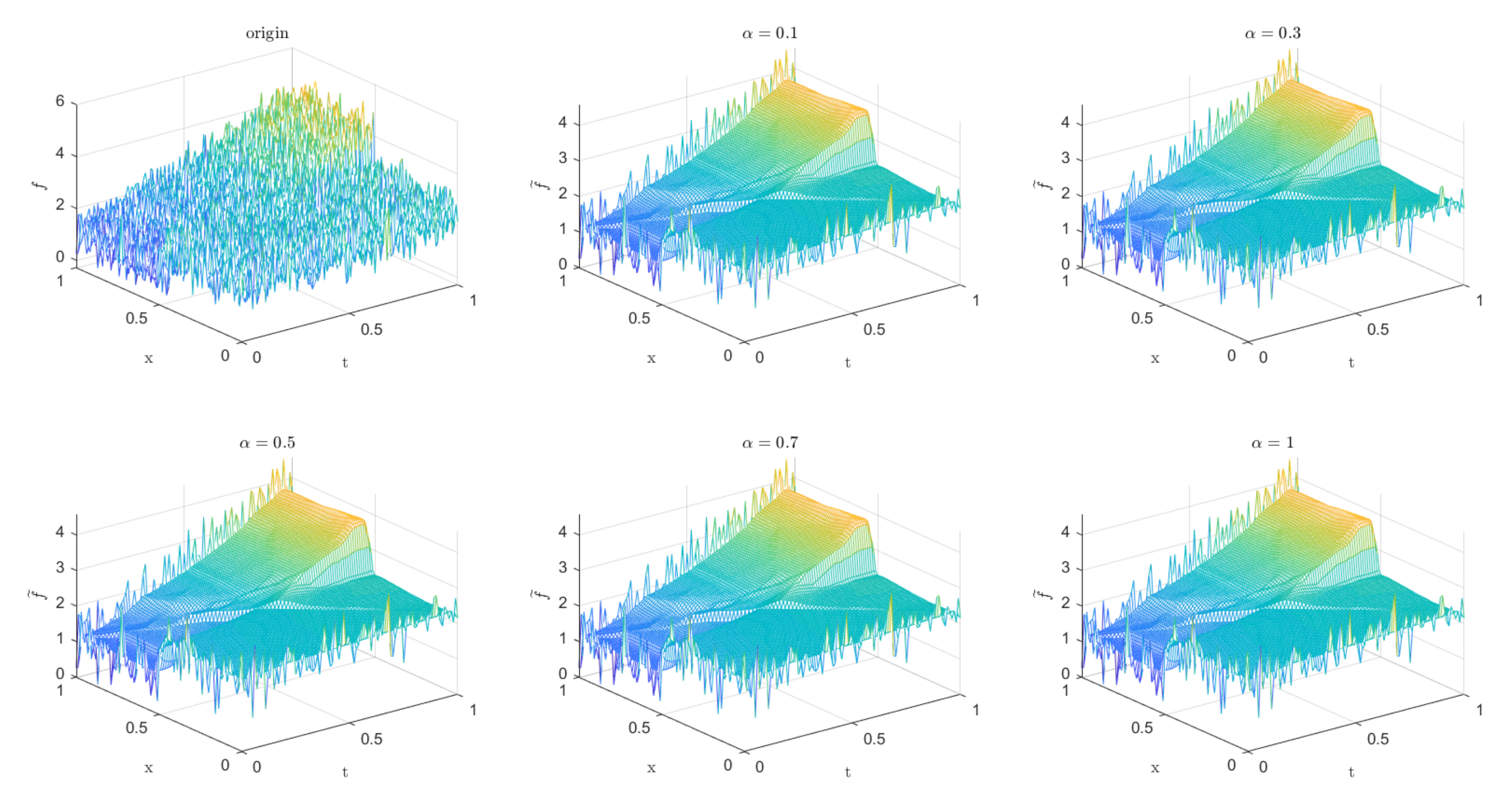}
	\caption{\textbf{Discover $\widetilde{f}$ for $h_{x}=h_{t}=1/100$ with $50\%$ Guassian noise in Type 1 for Example \ref{ex31}.  }}\label{6.2.1.1}
\end{figure*}

\begin{figure*}[!h]
	\centering
	\includegraphics[scale=0.35]{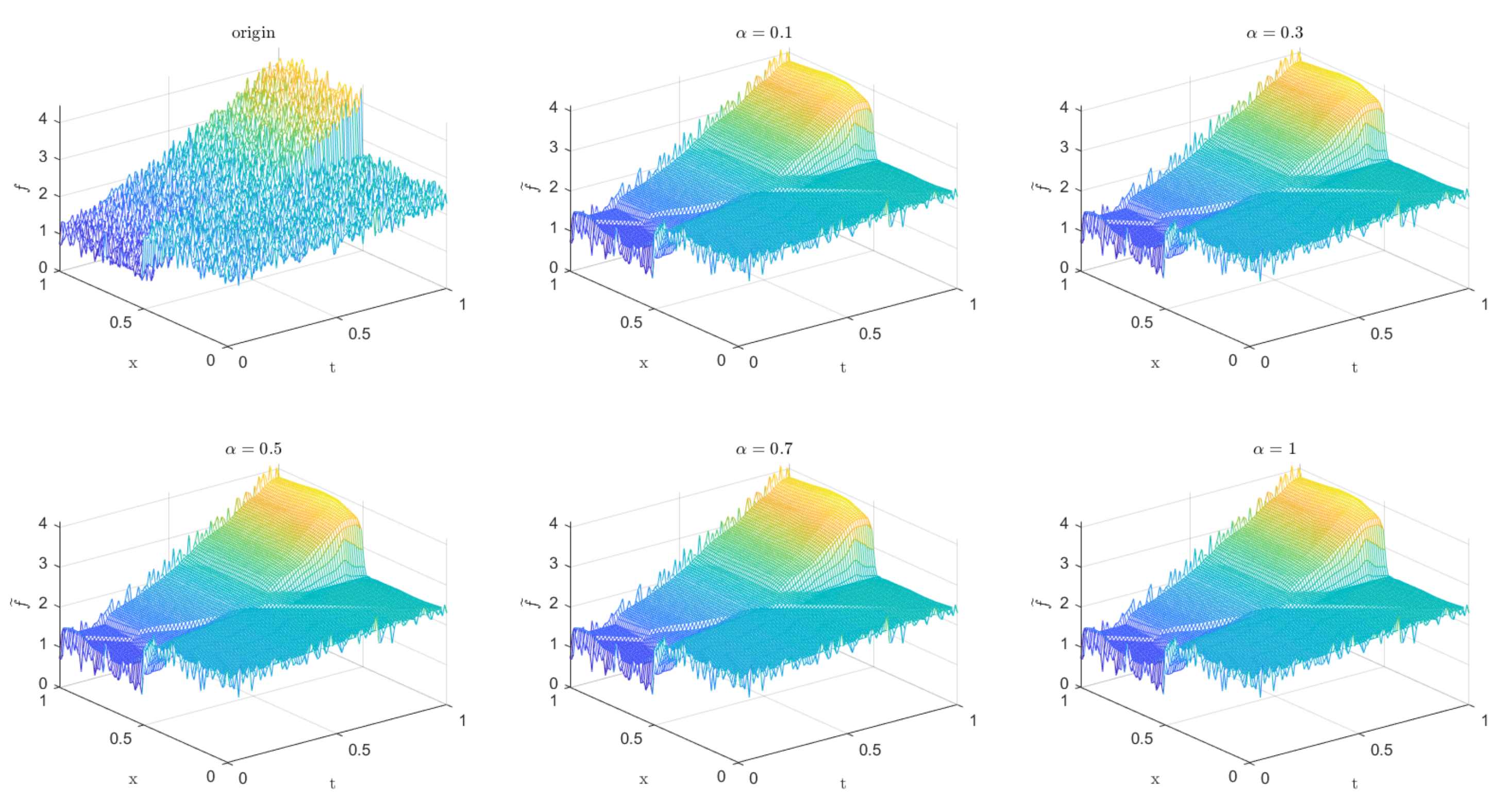}
	\caption{\textbf{Discover $\widetilde{f}$ for $h_{x}=h_{t}=1/100$ with $20\%$ Guassian noise in Type 1 for Example \ref{ex31}.  }}\label{6.2.1.2}
\end{figure*}

\begin{figure*}[!h]
	\centering
	\includegraphics[scale=0.35]{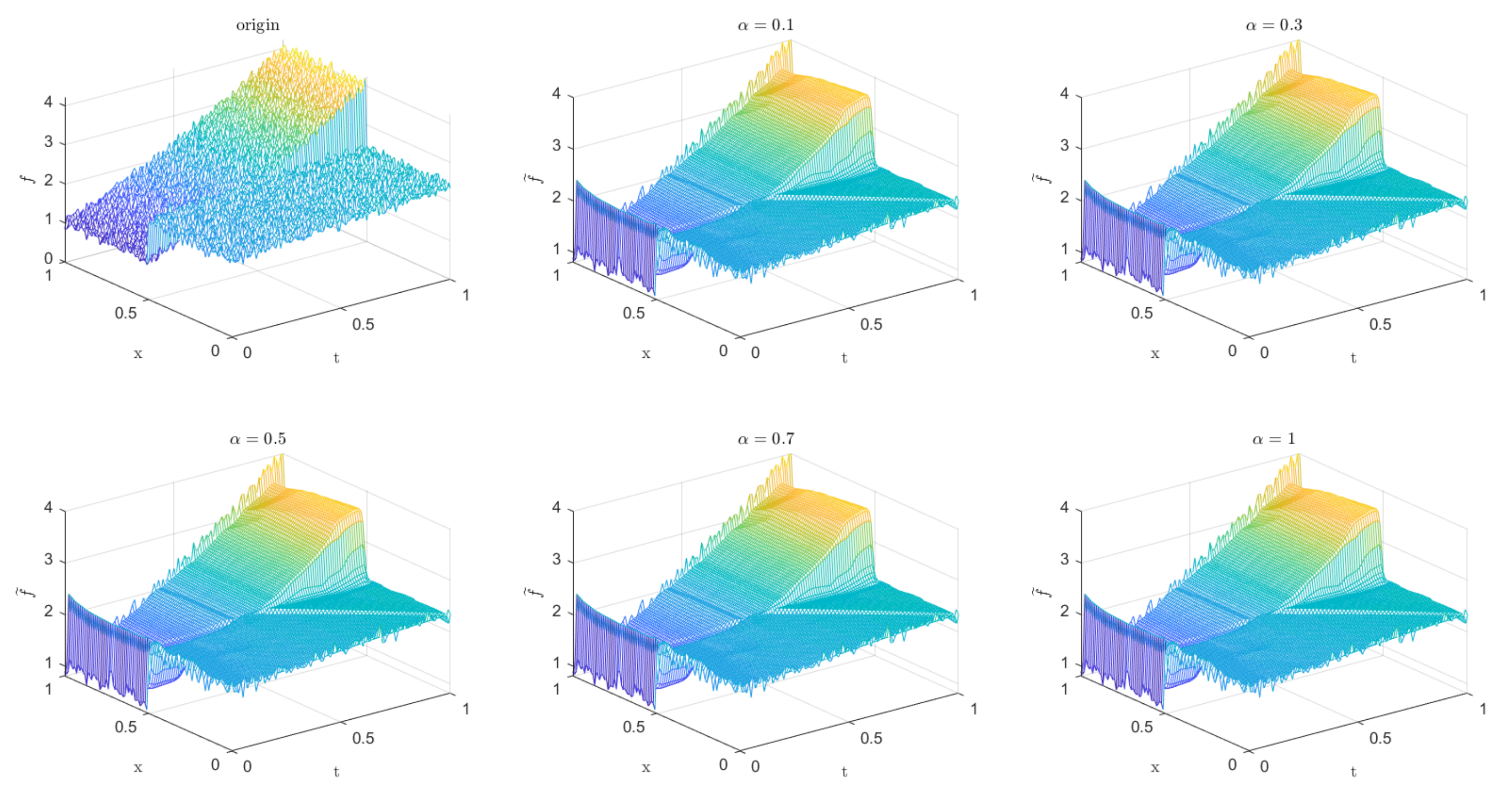}
	\caption{\textbf{Discover $\widetilde{f}$ for $h_{x}=h_{t}=1/100$ with $10\%$ Guassian noise in Type 1 for Example \ref{ex31}.  }}\label{6.2.1.3}
\end{figure*}

\begin{figure*}[!h]
	\centering
	\includegraphics[scale=0.35]{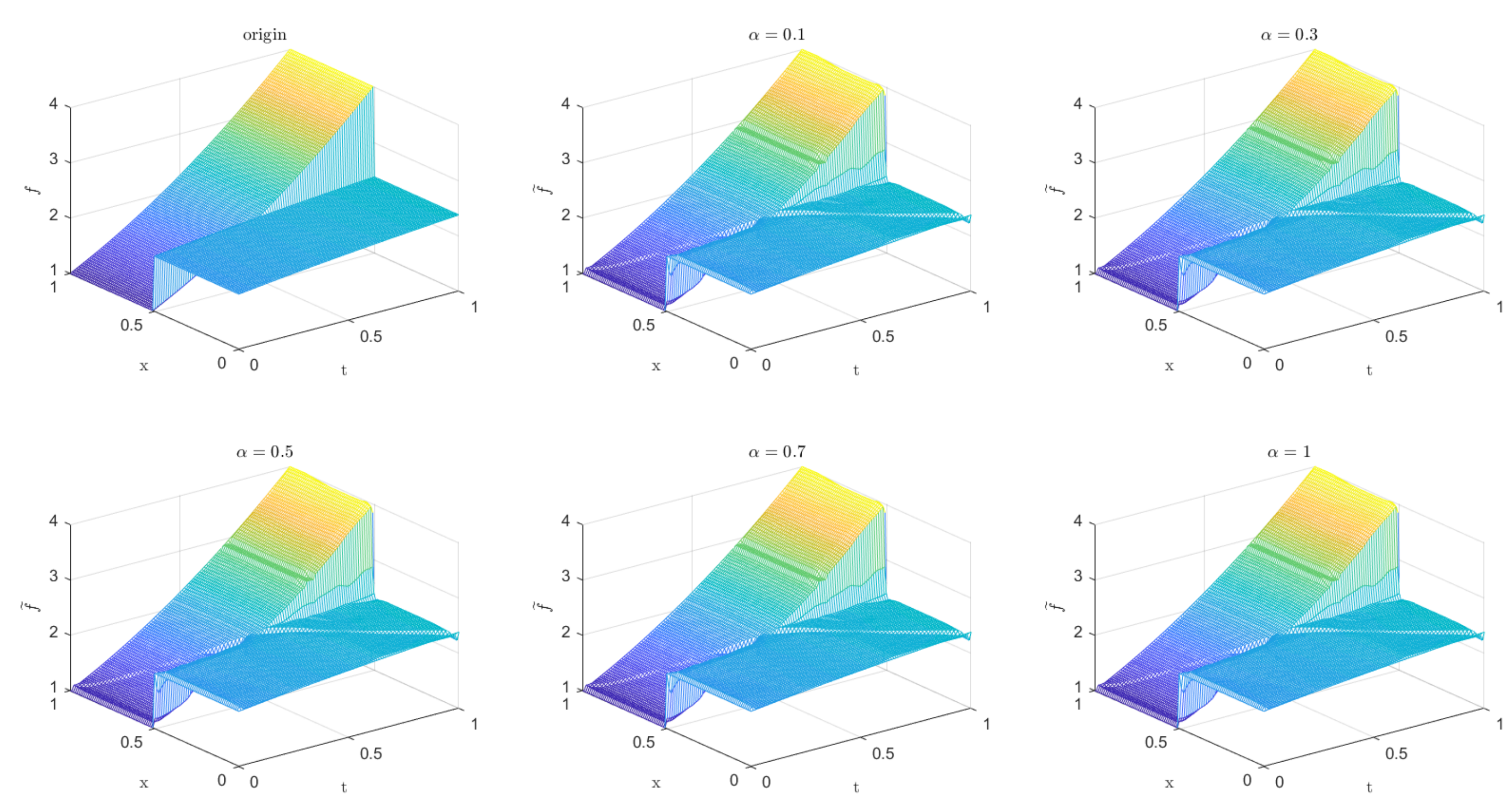}
	\caption{\textbf{Discover $\widetilde{f}$ for $h_{x}=h_{t}=1/100$ with clean data in Type 1 for Example \ref{ex31}.  }}\label{6.2.1.4}
\end{figure*}

Table  \ref{Tab5} and Figures \ref{6.2.1.1}-\ref{6.2.1.4} show the relative error \eqref{D6.1} between the source function $f$ in Example \ref{ex21} and discovery $\widetilde{f}$ in  \eqref{C6}
with different noise levels and epoch=$255$ in Type 1. The numerical experiments  are given to illustrate  the availability using deep learning for the different noise levels.
From Table \ref{Tab5}, it can be found that our proposed algorithm is stable and accurate faced with different discretization sizes and noise level,
which is robust  even for  $50\%$  Guassian noise level. From  Figures  \ref{6.2.1.1}-\ref{6.2.1.4}, it also appears a  layer or blows up at $t=0$ (boundary noise pollution).

\begin{table}[!ht]
	\centering
	\renewcommand\arraystretch{1.5}
	\caption{\textbf{The relative error  $||e_{\hat{f}}||_{r}$ with Guassian noise  in Type 2 for Example \ref{ex31}.   }} \label{Tab5}
	\label{tab:the error comparison}
	\resizebox{\textwidth}{!}{
	\begin{tabular}{|c|c|ccccc|}
		\hline
	Threshold& \diagbox{$h_{x}$}{$\alpha$} & 0.1   &0.3  &0.5 &0.7 &1\\
	\hline
		   	 \multirow{3}{*}{$50\%$ noise}&1/25 &2.1181e-01&	2.1181e-01&2.1181e-01&2.1181e-01&2.1181E-01\\
		\cline{2-7}
		&1/50 &2.1048e-01&2.1048e-01&2.1048e-01&2.1048e-01	&2.1048e-01\\
		\cline{2-7}
	   &1/100&2.1618e-01&2.1618e-01&2.1618e-01&2.1618e-01	&2.1618e-01\\
	   \hline
	   	   	 \multirow{3}{*}{$20\%$ noise}&1/25 &1.1250e-01&	1.1250e-01&1.1250e-01&1.1250e-01&1.1250e-01
\\
		\cline{2-7}
		&1/50 &9.5174e-02&9.5174e-02&9.5174e-02&9.5174e-02&	9.5174e-02\\
		\cline{2-7}
	   &1/100&9.1433e-02&9.1433e-02&9.1433e-02&9.1433e-02&	9.1433e-02\\
	   \hline
	    \multirow{3}{*}{$10\%$ noise}&1/25 &1.1917e-01&1.1917e-01	&1.1917e-01&1.1917e-01&1.1917e-01\\
		\cline{2-7}
		&1/50 &8.6597e-02&8.6597e-02&8.6597e-02&8.6597e-02&	8.6597e-02\\
		\cline{2-7}
	   &1/100&7.4511e-02&7.4511e-02&7.4511e-02&7.4511e-02&	7.4511e-02\\
	   \hline
	    \multirow{3}{*}{clean data}&1/25 &7.0042e-02&7.0042e-02&	7.0042e-02&7.0042e-02&7.0042e-02\\
		\cline{2-7}
		&1/50 &4.0945e-02&4.0945e-02&4.0945e-02&4.0945e-02&	4.0945e-02\\
		\cline{2-7}
	   &1/100&1.4316e-02&1.4316e-02&1.4316e-02&1.4316e-02&	1.4316e-02\\
	   \hline

	\end{tabular}}
\end{table}

\subsubsection{Type 2 in DNN}
\begin{figure*}[!h]
	\centering
	\includegraphics[scale=0.35]{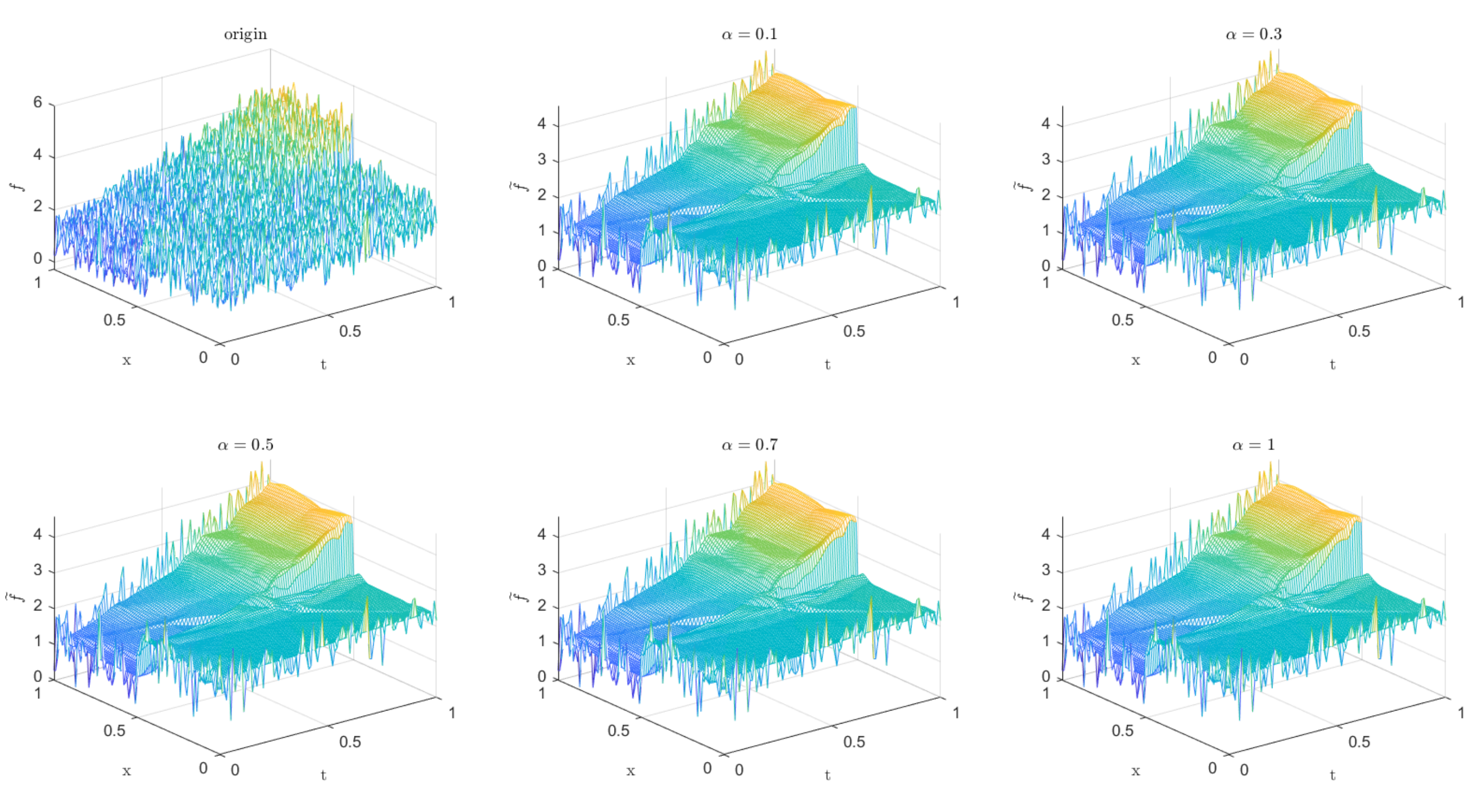}
	\caption{\textbf{Discover $\widetilde{f}$ for $h_{x}=h_{t}=1/100$ with $50\%$ Guassian noise in Type 2 for Example \ref{ex31}.  }}\label{6.2.2.1}
\end{figure*}

\begin{figure*}[!h]
	\centering
	\includegraphics[scale=0.35]{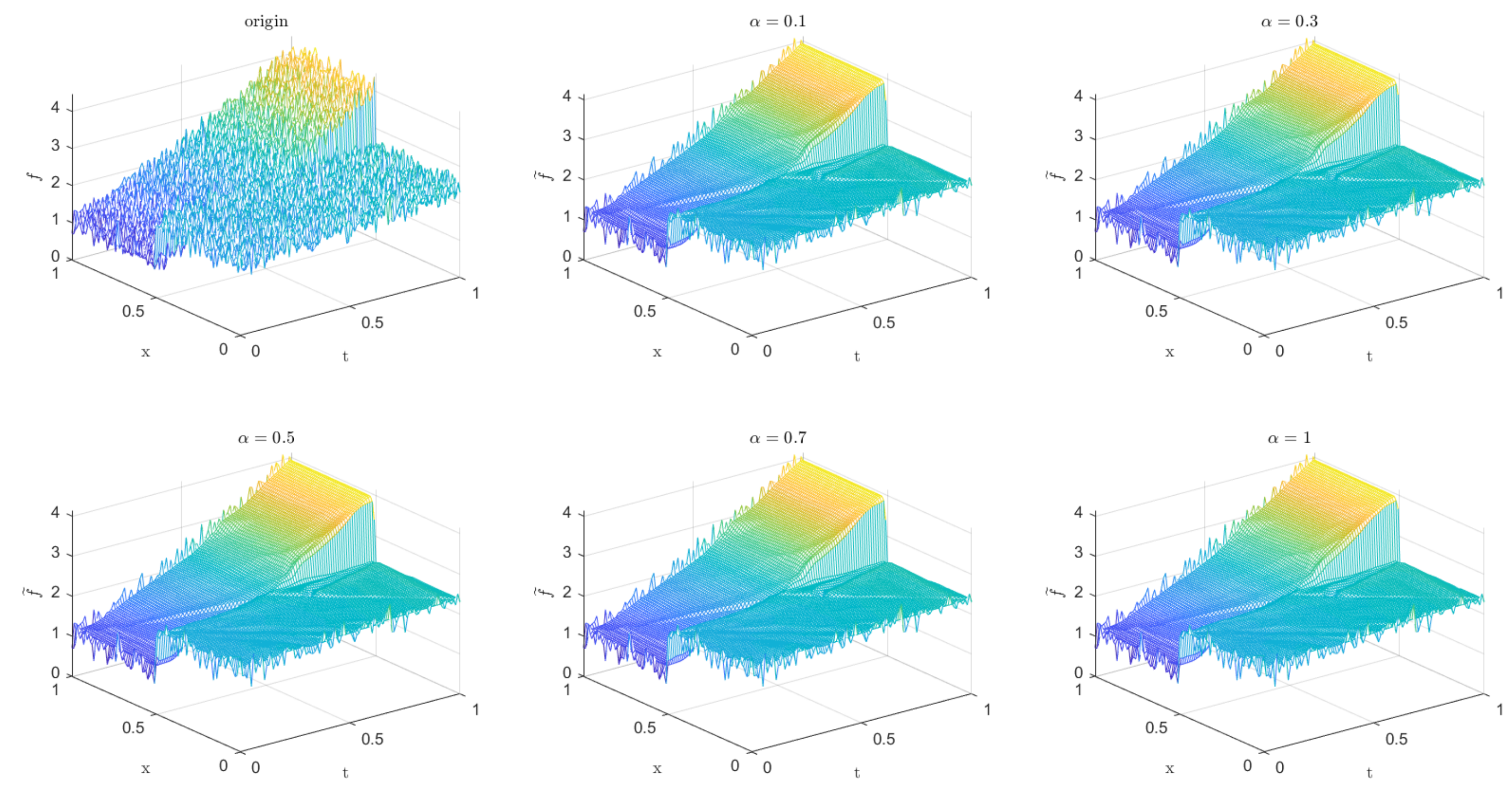}
	\caption{\textbf{Discover $\widetilde{f}$ for $h_{x}=h_{t}=1/100$ with $20\%$ Guassian noise in Type 2 for Example \ref{ex31}.  }}\label{6.2.2.2}
\end{figure*}

\begin{figure*}[!h]
	\centering
	\includegraphics[scale=0.35]{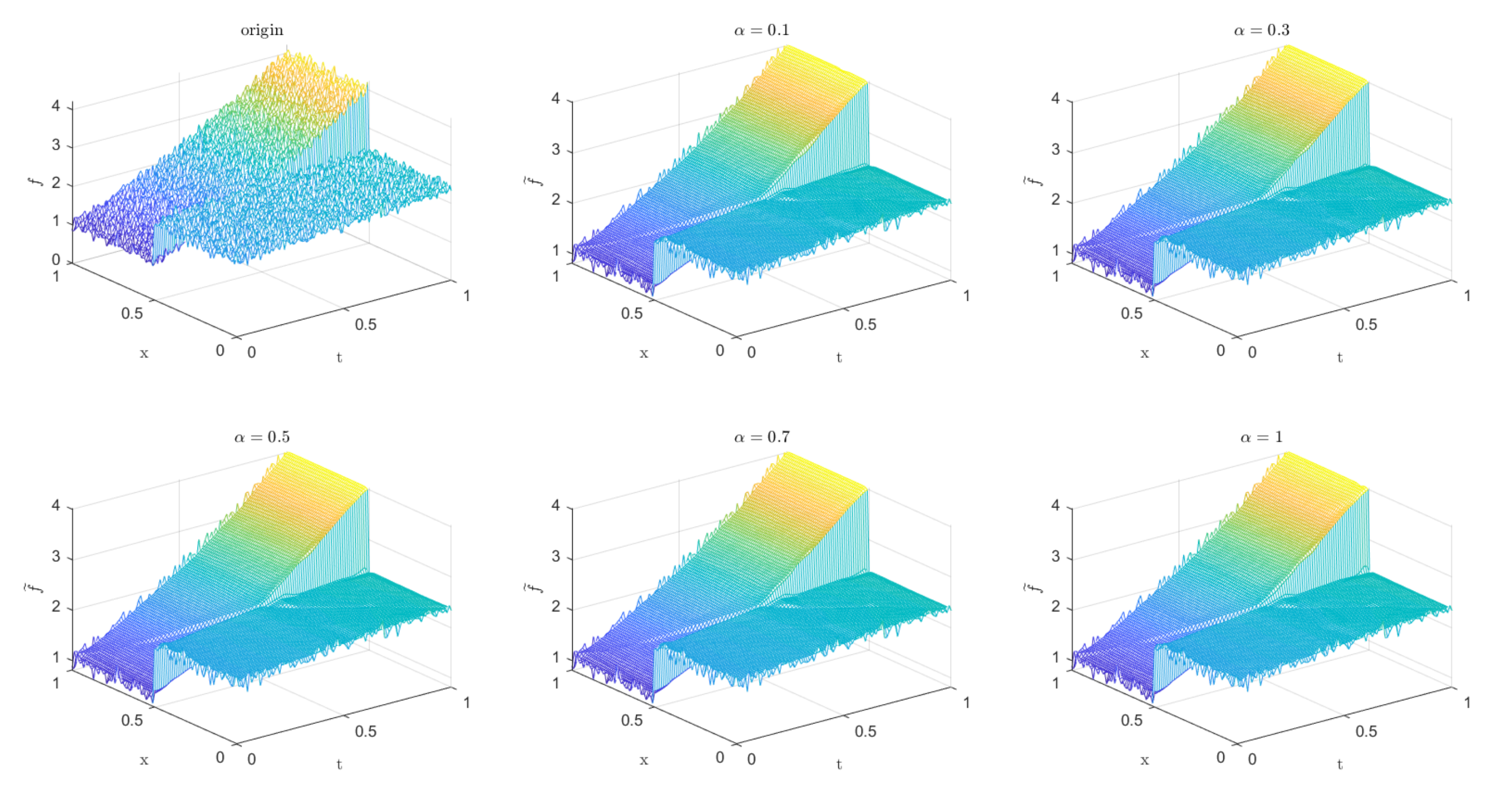}
	\caption{\textbf{Discover $\widetilde{f}$ for $h_{x}=h_{t}=1/100$ with clean data in Type 2 for Example \ref{ex31}.  }}\label{6.2.2.3}
\end{figure*}

Table  \ref{Tab6} and Figures \ref{6.2.2.1}-\ref{6.2.2.3} show the relative error \eqref{D6.1} between the source function $f$ in Example \ref{ex21} and discovery $\widetilde{f}$ in  \eqref{C6}
with different noise levels and epoch=$270$ in Type 2. The numerical experiments  are given to illustrate  the availability using deep learning for the different noise levels.
From Table \ref{Tab6}, it can be found that our proposed algorithm is stable and accurate faced with different discretization sizes and noise level,
which is robust  even for  $50\%$ Guassian noise level.  From  Figures  \ref{6.2.2.1}-\ref{6.2.2.3}, it also appears a littlt layer or blows up at $t=0$ (boundary noise pollution).
\begin{table}[!ht]
	\centering
	\renewcommand\arraystretch{1.5}
	\caption{\textbf{The relative error  $||e_{\hat{f}}||_{r}$ with Guassian noise  in Type 2 for Example \ref{ex31}.   }} \label{Tab6}
	\label{tab:the error comparison}
	\resizebox{\textwidth}{!}{
	\begin{tabular}{|c|c|ccccc|}
		\hline
	Threshold& \diagbox{$h_{x}$}{$\alpha$} & 0.1   &0.3  &0.5 &0.7 &1\\
	\hline
		   	 \multirow{3}{*}{\tabincell{c}{$h^{*}_{x}=h^{*}_{t}$\\$=1/100$ with\\$50\%$ noise}}&1/25 &2.1509e-01&	2.1435e-01&2.1464e-01&2.1482e-01&2.1413e-01\\
		\cline{2-7}
		&1/50 &2.1005e-01&2.0957e-01&2.0972e-01&2.0975e-01	&2.0946e-01\\
		\cline{2-7}
	   &1/100&2.0500e-01&2.0510e-01&2.0529e-01&2.0529e-01&	2.0537e-01\\
	   \hline
	   	   	 \multirow{3}{*}{\tabincell{c}{$h^{*}_{x}=h^{*}_{t}$\\$=1/100$ with\\$20\%$ noise}}&1/25 &1.0418e-01	&1.0426e-01&1.0441e-01&1.0465e-01&1.0478e-01\\
		\cline{2-7}
		&1/50 &9.2908e-02&9.2833e-02&9.2892e-02&9.3044e-02	&9.2955e-02\\
		\cline{2-7}
	   &1/100&8.5171e-02&8.4946e-02&8.5225e-02&8.5757e-02&	8.8259e-02\\
	   \hline
	    \multirow{3}{*}{\tabincell{c}{$h^{*}_{x}=h^{*}_{t}$\\$=1/100$ with\\$10\%$ noise}}&1/25 &7.3662e-02&7.4271e-02&	7.4214e-02&7.4147e-02&7.4018e-02\\
		\cline{2-7}
		&1/50 &5.7890e-02&5.8266e-02&5.8156e-02&5.8111e-02	&5.8036e-02\\
		\cline{2-7}
	   &1/100&4.2822e-02&4.2814e-02&4.2745e-02&4.2785e-02	&4.2809e-02\\
	   \hline
	    \multirow{3}{*}{\tabincell{c}{$h^{*}_{x}=h^{*}_{t}$\\$=1/100$ with\\$0\%$ noise}}&1/25 &5.6927e-02&5.8353e-02	&5.8272e-02&5.7993e-02&5.8502e-02\\
		\cline{2-7}
		&1/50 &3.8296e-02&3.9283e-02&3.9420e-02&3.9318e-02&	3.9337e-02\\
		\cline{2-7}
	   &1/100&5.4745e-03&4.4494e-03&4.9254e-03&4.7000e-03	&3.9269e-03\\
	   \hline
	\end{tabular}}
\end{table}

\section{Conclusion}
In this work, we firstly propose the discovery of subdiffusion with noisy data  in deep learning based on the designing  two types.
The numerical experiments  are given to illustrate  the availability using deep learning even for the $50\%$ noise level.
The advantage of  first Type  is that model \eqref{C6} can be  trained    for each fixed $\alpha$, which reduces the computational count and required storage, since it fixes  $\alpha$ as input data.
Obviously, it may loss some accuracy. To recover the accuracy,  second Type  is complemented, which  needs more computational count.
It implies  an interesting generalized structure by combining  Type 1 and Type 2, which may keep suitable accuracy and reduces computational count.
As a result, even when discovering the interger-order equations, the proposed type DNNs have some advantages compared with the traditional numerical scheme and previous works to solve the data-driven discovery for differential equations.
The interesting topic is reducing the computational count and storage by fast multigrid or  conjugate gradient squared method \cite{CCNWL:20,CES:20}.
Another interesting topic is how to design the correction of  $L_1$ approximation for
reducing the  boundary noise pollution, since subdiffuison model appears a layer or blows up at $t=0$ \cite{ShCh:2020,SOG:17} with low time regularity.

\bibliographystyle{amsplain}

\begin{thebibliography}{10}

\bibitem{CCNWL:20}
Cao,  R.J.,   Chen, M.H.,  Ng,  M.K.,    Wu, Y.J.:
Fast and high-order accuracy numerical methods for  time-dependent nonlocal problems in $\mathbb{R}^2$. 
 J. Sci. Comput. \textbf{84:8} (2020).




\bibitem{CJR:2021} Chen, J.R., Jin, S., Lyu, L.Y.:
A Deep learning Based Discotinuous Galerkin Method for Hyperbolic Equations with Discontinuous Solutions and Random Uncertainties. arXiv:2107.01127




\bibitem{ChJB:21}    Chen, M.H., Jiang, S.Zh., Bu, W.P.:
Two $L1$ schemes on graded meshes for fractional Feynman-Kac equation.
J. Sci. Comput. \textbf{88:58} (2021)


\bibitem{CES:20}
Chen, M.H., Ekstr\"{o}m, S.E., Serra-Capizzano, S.:
 A Multigrid method for nonlocal problems: non-diagonally dominant or Toeplitz-plus-tridiagonal systems.
  SIAM J. Matrix Anal. Appl. \textbf{41} 1546-1570 (2020)
\bibitem{CWHZ:21} Chen, W.Q., Wang, Q., Hesthaven, J.S., Zhang, C.H.:
Physics-informed machine learning for reduced-order modeling of nonlinear problems. J. Comput. Phys. \textbf{446}, 110666 (2021)

\bibitem{EK:11} Eliazar, I., Klafter, J.: Anomalous is ubiquitous.
Ann. Physies \textbf{326}, 2517--2531 (2011)



\bibitem{DUQ:2021} Du. Q., Gu, Y.Q., Yang, H.Z., Zhou, C.: The discovery of dynamics via linear multistep methods and deep learning:error estimation. arXiv:2103.11488



\bibitem{CGD:2021}Duan, C.G., Jiao, Y.L, Lai, Y.M., Lu, X.L.,Yang, Z.J.: Convergence rate analysis for deep Ritz method. arXiv:2103.13330


\bibitem{WEBY:2017} E, W.N., Yu, B.: The Deep Ritz method:A deep learning-based numerical algorithm for solving variational problems. Commun. Math. Stat. \textbf{6}, 1--12 (2018)




\bibitem{IYA:2016} Goodfellow, I., Bengio, Y., Courville, A.:
Deep Learning. Adaptive Computation and Machine Learning.  MIT Press, Cambridge (2016)

\bibitem{YGMK:2021} Gu, Y.Q., Ng, M.K.: Deep Ritz method for the spectral fractional Laplacian equation using the Caffarelli-Silvestre extension. arXiv:2108.11592


\bibitem{GRPK:19} Gulian, M.K., Raissi, M., Perdikaris, P., Karniadakis, G.:
Machine learning of space-fractional differential equations.
SIAM J. Sci. Comput. \textbf{41}, A2485--A2509 (2019)


\bibitem{YLJ:2021} Jiao, Y.L., Lai, Y.M., Lo,Y.S., Wang.Y., Yang.Y.F.: Error Analysis of Deep Ritz Methods for Elliptic Equations. arXiv:2107.14478




\bibitem{RD:21} Keller, R., Du, Q.: Discovery of dynamics using linear multistep methods. SIAM J. Numer. Anal. \textbf{59}, 429--455 (2021)


\bibitem{LWD:17}  Li, Y.J., Wang, Y.J., Deng, W.H.:
Galerkin finite element approximations for stochastic space-time fractional wave equations.
SIAM J. Numer. Anal.  \textbf{55},  3173--3202 (2017)

\bibitem{LZY:2021} Li, Z.Y., Kovachki, N., Azizzadenesheli., K., Liu, B. ,Bhattacharya, K., Stuart, A., Anandkumar, A.:
Markov Neurl Operators for Learning Chaotic Systems. arXiv:2106.06898


\bibitem{Lin:07}  Lin, Y.M., Xu, C.J.:
Finite difference/spectral approximations for the time-fractional diffusion equation.
J. Comput. Phys. \textbf{225}, 1533--1552 (2007)


\bibitem{MS:12} Meerschaert, M.M., Sikorskii, A.:
Stochastic Models for Fractional Calculus. de Gruyter GmbH Berlin  (2012)

\bibitem{M2000} Metzler R., Klafter J.: The random walk's guide to anomalous diffusion: a fractional dynamics approach. Phys. Rep. \textbf{339}, 1--77 (2000)







\bibitem{Podlubny:99} Podlubny I.:  Fractional Differential Equations. Academic Press, New York (1999)


\bibitem{QWX:2019} Qin, T., Wu, K., Xiu, D.: Data driven governing equations approximation using deep neural network. J. Comput. Phys. \textbf{395}, 620--635 (2019)




\bibitem{MRI:2018} Raissi, M.: Deep hidden physics models deep learning of nonlinear partial differential equations. The Journal of Machine Learning Research. \textbf{19}, 1--24 (2018)


\bibitem{SHR:2019}  Rudy, S.H., Kutz, J.N, Brunton, S.L.: Deep learning of dynamics and signal-noise decompoisition with time-stepping constraints. J. Comput. Phys. \textbf{396}, 483--506 (2019)



\bibitem{ZuowShen:2021} Shen, Z.W., Yang, H.Z., Zhang, S.J.:
Deep Network Approximation Characterized by Number of Neurons Commun. Comput. Phys., \textbf{28}, 1768--1811 (2019)


\bibitem{DEM:2020}
Shen, X., Cheng, X.L., Liang, .K.W.: Deep Euler Method: Solving ODEs by approximating the local truncation error of the euler method.  arXiv:2003.09573


\bibitem{DGM:2018} Sirignano, J., Spiliopoulos, K.:
DGM: A deep learning algorithm for solving partial differental equations. J. Comput. Phys. \textbf{375}, 1339--1364 (2018)

\bibitem{ShCh:2020}    {Shi, J.K., Chen, M.H.:}
{Correction of high-order BDF convolution quadrature for  fractional Feynman-Kac equation with L\'{e}vy flight}.
J. Sci. Comput. \textbf{85:28} (2020)


\bibitem{SOG:17}  Stynes, M., O'riordan, E., Gracia, J.L.:
Error analysis of a finite difference method on graded meshes for a time-fractional diffusion equation.
SIAM J. Numer. Anal.  \textbf{55},  1057--1079 (2017)



\bibitem{RTP:2019} Tipireddy, R., Perdikaris, P., Stinis , P., Tartakovsky, A.: A comparative study of physics-informed neural network models for learning unknown dynamics and constitutive relations. arXiv:1904.04058





\bibitem{WCDBD:22} Wang, C., Chen, M.H., Deng, W.H., Bu, W.P., Dai, X.J.:
A sharp error estimate of Euler--Maruyama method for stochastic Volterra integral equations.
Math. Method Appl. Sci.  Minor Revised.



\bibitem{XUY:2020} Xu, Y., Zhang, H., Li, Y.G.,  Zhou, K., Liu, Q., Jurgen, K.: Solving Fokker-Planck equation using deep learning. Chaos. \textbf{30}, 013133 (2020)







\bibitem{ZK:17} Zhang, Z.H., Karniadakis, G.E.:
Numerical Methods for Stochastic Partial Differential Equations with White Noise. Springer, New York (2017)







\end{thebibliography}

\end{document}